\documentclass[11pt]{article}

\usepackage{amssymb,amsbsy,amsmath,amscd,amsthm,amsfonts}
\usepackage{cite}
\usepackage{enumerate}
\usepackage{dsfont}
\usepackage{graphicx}
\usepackage{longtable}
\usepackage{tablefootnote}
\usepackage{standalone}
\usepackage{fancyhdr}
\usepackage{float}
\usepackage{caption}
\usepackage{subcaption}

\usepackage[margin=1in]{geometry}

\usepackage[
            CJKbookmarks=true,
            bookmarksnumbered=true,
            bookmarksopen=true,
            colorlinks=true,
            citecolor=red,
            linkcolor=blue,
            anchorcolor=red,
            urlcolor=blue
            ]{hyperref}

\newtheorem{theorem}{Theorem}

\newtheorem{lemma}{Lemma}

\newtheorem{assumption}{Assumption}

\usepackage{prettyref}
\newrefformat{eq}{(\ref{#1})}
\newrefformat{thm}{Theorem~\ref{#1}}
\newrefformat{sec}{Section~\ref{#1}}
\newrefformat{algo}{Algorithm~\ref{#1}}
\newrefformat{fig}{Fig.~\ref{#1}}
\newrefformat{tab}{Table~\ref{#1}}
\newrefformat{rmk}{Remark~\ref{#1}}
\newrefformat{def}{Definition~\ref{#1}}
\newrefformat{cor}{Corollary~\ref{#1}}
\newrefformat{lmm}{Lemma~\ref{#1}}
\newrefformat{prop}{Proposition~\ref{#1}}
\newrefformat{app}{Appendix~\ref{#1}}
\newrefformat{ex}{Example~\ref{#1}}
\newrefformat{ass}{Assumption~\ref{#1}}

\title{Risk Bounds for High-dimensional Ridge Function Combinations Including Neural Networks}
\date{}
\author{Jason M. Klusowski \qquad Andrew R. Barron \\ \\
Department of Statistics \\ Yale University \\ New Haven, CT, USA \\ \\
Email: \{jason.klusowski, andrew.barron\}@yale.edu}

\begin{document}

\maketitle

\tableofcontents

\abstract

Let $ f^{\star} $ be a function on $ \mathbb{R}^d $ with an assumption of a spectral norm $ v_{f^{\star}} $. For various noise settings, we show that $ \mathbb{E}\|\hat{f} - f^{\star} \|^2 \leq \left(v^4_{f^{\star}}\frac{\log d}{n}\right)^{1/3} $, where $ n $ is the sample size and $ \hat{f} $ is either a penalized least squares estimator or a greedily obtained version of such using linear combinations of sinusoidal, sigmoidal, ramp, ramp-squared or other smooth ridge functions. The candidate fits may be chosen from a continuum of functions, thus avoiding the rigidity of discretizations of the parameter space. On the other hand, if the candidate fits are chosen from a discretization, we show that $ \mathbb{E}\|\hat{f} - f^{\star} \|^2 \leq \left(v^3_{f^{\star}}\frac{\log d}{n}\right)^{2/5} $.

This work bridges non-linear and non-parametric function estimation and includes single-hidden layer nets. Unlike past theory for such settings, our bound shows that the risk is small even when the input dimension $ d $ of an infinite-dimensional parameterized dictionary is much larger than the available sample size. When the dimension is larger than the cube root of the sample size, this quantity is seen to improve the more familiar risk bound of $ v_{f^{\star}}\left(\frac{d\log (n/d)}{n}\right)^{1/2} $, also investigated here.

\section{Introduction}

Functions $ f^{\star} $ in $ \mathbb{R}^d $ are approximated using linear combinations of ridge functions with one layer of nonlinearities. These approximations are employed via functions of the form
\begin{equation} \label{eq:general}
f_m(x) = f_m(x,\zeta) = \sum_{k=1}^m c_k \phi(a_k \cdot x+b_k),
\end{equation}
which is parameterized by the vector $ \zeta $, consisting of $ a_k $ in $ \mathbb{R}^d $, and $ b_k, c_k $ in $ \mathbb{R} $ for $ k = 1,\dots,m $, where $ m \geq 1 $ is the number of nonlinear terms. Models of this type arise with considerable freedom in the choice of the activation function $ \phi $, ranging from general smooth functions of projection pursuit regression \cite{Friedman1981} to the unit step sigmoid and ramp functions of single-hidden layer neural nets \cite{Barron1994, Barron1992, Barron1993, Breiman1993, Lee1996}.

Our focus in this paper is on the case that $ \phi $ is a fixed Lipschitz function (such as a sigmoid or ramp or sinusoidal function), though some of our conclusions apply more generally. For these activation functions, we will obtain statistical risk conclusions using a penalized least squares criterion. We obtain generalization error bounds for these by balancing the approximation error and descriptive complexity. The most general form of our bounds hold for quite general non-linear infinite dictionaries. A hallmark of our conclusions is to lay bare how favorable risk behavior can be obtained as long as the logarithm of the number of parameters relative to sample size is small. This entails a slower rate of convergence through a rate that is smaller than what is cemented in traditional cases, but leads to better results than these earlier bounds would permit in certain very high-dimensional situations. From an applied perspective, good empirical performance of neural net (and neural net like) models has been reported as in \cite{LeCun2015} even when $ d $ is much larger than $ n $, though theoretical understanding has been lacking. Returning to the case of a single layer of nonlinearly parameterized function, it is useful to view the representation \prettyref{eq:general} as
\begin{equation} \label{eq:simple_general}
\sum_{h} \beta_h h(x),
\end{equation}
where the $ h $ are a selection of functions from the infinite library $ \mathcal{H} = \mathcal{H}_{\phi} $ of functions of the form $ \pm \phi(\theta \cdot x) $ for real vector $ \theta $ and the $ \beta_h $ are coefficients of linear combination of $ \pm \phi $ in the library. These representations are single hidden-layer networks. Deep network approximations are not very well understood. Nevertheless our results generalize provided some of our arguments are slightly modified.

We can reduce \prettyref{eq:general} to \prettyref{eq:simple_general} as follows. Suppose the library is symmetric $ \mathcal{H} = -\mathcal{H} $ and contains the zero function. Without loss of generality, we may assume that the $ c_k $ or $ \beta_h $ are non-negative by replacing the associated $ \phi $ with $ \phi\ \text{sgn}c_k $, that by assumption also belongs to $ \mathcal{H} $. One can assume the internal parameterization $ a \cdot x+b $ take the form $ \theta \cdot x $ by appending a coordinate of constant value $ 1 $ to $ x $ and a coordinate of value $ b $ to the vector $ a $. Note that now $ x $ and $ \theta $ are $ (d+1) $-dimensional.

We will take advantage of smoothness of activation function (assumption that either $ \phi $ is Lipschitz or that its first derivative $ \phi' $ is Lipschitz). Suppose $ P $ is an arbitrary probability measure on $ [-1,1]^d $. Let $ \|\cdot\| $ be the $ L^2(P) $ norm induced by the inner product $ \langle \cdot, \cdot \rangle $. 
For a symmetric collection of dictionary elements $ \mathcal{H} = -\mathcal{H} $ containing the zero function, we let $ \mathcal{F} = \mathcal{F}_{\mathcal{H}} $ be the linear span of $ \mathcal{H} $.

The variation $ v_f = \|f\|_{\mathcal{H}} $ of $ f $ with respect to $ \mathcal{H} $ (or the atomic norm of $ f $ with respect to $ \mathcal{H} $) is defined by
\begin{equation*}
\lim_{\delta\downarrow 0}\inf_{f_{\delta}\in\mathcal{F}}\left\{  \|\beta\|_1 : f_{\delta} = \sum_{h\in\mathcal{H}}\beta_{h}h \; \text{and} \; \|f_{\delta} - f \| \leq \delta,\; \beta_{h}\in \mathbb{R}^{+} \right\},
\end{equation*}
where $ \|\beta\|_1 = \sum_{h\in\mathcal{H}}\beta_h $. For functions in $ \mathcal{F}_{\mathcal{H}} $, this variation picks out the smallest $ \|\beta\|_1 $ among representations $ f = \sum_{h\in\mathcal{H}}\beta_{h}h $. In the particular case that $ f = \sum_{h\in\mathcal{H}}\beta_h $, we have $ v_f = \|\beta\|_1 $.  For functions in the $ L^2(P) $ closure of the linear span of $ \mathcal{H} $, the variation is the smallest limit of such $ \ell_1 $ norms among functions approaching the target. The subspace of functions with $ \|f\|_{\mathcal{H}} $ finite is denoted $ L_{1,\mathcal{H}} $. Such variation control provides for approximation (opportunity) for dimension independent rates of order $ 1/\sqrt{m} $ with an $ m $ term approximation.

It is fruitful to discuss spectral conditions for finite variation for various choices of $ \phi $. To this end, define $ v_{f^{\star},s} = \int_{\mathbb{R}^d}\|\omega\|^s_1\widetilde{f}(\omega)d\omega $, for $ s \geq 0 $. If $ f^{\star} $ has a bounded domain in $ [-1,1]^d $ and a Fourier representation $ f^{\star}(x) = \int_{\mathbb{R}^d}e^{i\omega\cdot x}\widetilde{f}(\omega)d\omega $ with $ v_{f^{\star},1} < +\infty $, it is possible to use approximating functions of the form \prettyref{eq:general} with a single activation function $ \phi $. Such activation functions $ \phi $ can be be general bounded monotone functions. We use $ x $ for vectors in $ \mathbb{R}^d $ and $ z $ for scalars such as $ z = \theta \cdot x $. As we have said, to obtain risk bounds in later sections, we will assume that either $ \phi $ is bounded Lipschitz or that, additionally, its derivative $ \phi' $ is Lipschitz. These two assumptions are made precise in the following statements.

\begin{assumption} \label{ass:ass1}
The activation function $ \phi $ has $ L_{\infty} $ norm at most one and satisfies
\begin{equation*} \label{eq:Lipschitz}
|\phi(z) - \phi(\tilde{z})| \leq L_1|z-\tilde{z}|,
\end{equation*}
for all $ z, \tilde{z} $ in $ \mathbb{R} $ and for some positive constant $ L_1 > 0 $. 
\end{assumption}

\begin{assumption} \label{ass:ass2}
The activation function $ \phi $ has $ L_{\infty} $ norm at most one and satisfies
\begin{equation*} \label{eq:Lipschitz}
|\phi(z) - \phi(\tilde{z})| \leq L_1|z-\tilde{z}|,
\end{equation*}
and
\begin{equation*} \label{eq:Smooth}
|\phi'(z) - \phi'(\tilde{z})| \leq L_2|z-\tilde{z}|,
\end{equation*}
for all $ z, \tilde{z} $ in $ \mathbb{R} $ and for some positive constants $ L_1 > 0 $ and $ L_2 > 0 $.
\end{assumption}

In particular, \prettyref{ass:ass2} implies that
\begin{equation*}
|\phi(z) - \phi(\tilde{z}) - (z-\tilde{z})\phi'(\tilde{z})| \leq \frac{1}{2}(z-\tilde{z})^2L_2,
\end{equation*}
for all $ z, \tilde{z} $ in $ \mathbb{R} $.

A result from \cite{Barron1993} provides a useful starting point for approximating general functions $ f^{\star} $ by linear combinations of such objects. Suppose $ v_{f^{\star},1} $ is finite. Then by \cite{Barron1993} the function $ f^{\star} $ has finite variation with respect to step functions and, consequently, there exists an artificial neural network of the form \prettyref{eq:general} with $ \phi(x) = \text{sgn}(x) $, $ \|a_k\|_1 = 1 $, and $ |b_k| \leq 1 $ such that, if a suitable constant correction is subtracted from $ f^{\star} $, then
\begin{equation*}
\|f^{\star}-f_m\|^2 \leq \frac{v^2_{f^{\star},1}}{m}.
\end{equation*}
In particular, $ f^{\star} $ minus a constant correction has variation less than $ v_{f^{\star},1} $.

If $ \phi $ has right at left limits $ -1 $ and $ +1 $, respectively, the fact that $ \phi(\tau x) \rightarrow \text{sgn}(x) $ as $ \tau \rightarrow +\infty $ allows one to use somewhat arbitrary activation functions as basis elements. For our results, it in undesirable to have unbounded weights. Accordingly, it is natural to impose a restriction on the size of the internal parameters and to also enjoy a certain degree of smoothness not offered by step functions. Although, it should be mentioned that classical empirical process theory allows one to obtain covering numbers for indicators of half-spaces (which are scale invariant in the size of the weights) by taking advantage of their combinatorial structure \cite{Bartlett2009}. Nevertheless, we adopt the more modern approach of working with smoothly parameterized dictionaries. In this direction, we consider the result in \cite{Breiman1993}, which allows one to approximate $ f^{\star} $ by linear combinations of ramp ridge functions (also known as first order ridge splines or hinging hyper-planes) $ (x\cdot \alpha-t)_{+} = \max\{0, x\cdot \alpha-t \} $, with $ \|\alpha\|_1 = 1 $, $ |t| \leq 1 $. 

The ramp activation function $ \phi(x) = (x)_{+} $ (also called a lower-rectified linear unit or ReLU) is currently one of the most popular form of artificial neural network activation functions, particularly because it is continuous and Lipschitz. In particular, it satisfies the conditions of \prettyref{ass:ass1} with $ L_1 = 1 $ depending on the size of its domain.
In \prettyref{thm:splineapprox}, we refine a result from \cite{Breiman1993}. For an arbitrary target function $ f^{\star} $ with $ v_{f^{\star},2} $ finite has finite variation with respect to the ramp functions and, consequently, there exists an approximation of the form \prettyref{eq:general} activated by ridge ramp functions with $ \|a_k\|=1 $ and $ |b_k| \leq 1 $ such that if a suitable linear correction is subtracted from $ f^{\star} $, then
\begin{equation}
\|f^{\star}-f_m\|^2 \leq \frac{16v^2_{f^{\star},2}}{m}.
\end{equation}
In particular, $ f^{\star} $ minus a linear correction has variation less than $ v_{f^{\star},2} $.

The second order spline $ \phi(x) = (x)^2_{+} $, which may also be called ramp-squared, satisfies the conditions of \prettyref{ass:ass2} with constants $ L_1 = 2 $ and $ L_2 = 2 $ depending on the size of its domain. Likewise, in \prettyref{thm:splineapprox}, we show that for an arbitrary target function $ f^{\star} $ with $ v_{f^{\star},3} $ finite a quadratically corrected $ f^{\star} $ has finite variation with respect to the second order splines, there exists an approximation of the form \prettyref{eq:general} activated by second order ridge splines with $ \|a_k\|=1 $ and $ |b_k| \leq 1 $ such that, if a suitable quadratic correction is subtracted from $ f^{\star} $, then
\begin{equation}
\|f^{\star}-f_m\|^2 \leq \frac{16v^2_{f^{\star},3}}{m}.
\end{equation}
In particular, $ f^{\star} $ minus a quadratic correction has variation less than $ v_{f^{\star},3} $.

For integer $ s \geq 1 $, we define the infinite dictionary
\begin{align*}
\mathcal{H}_{s} = \{ x\mapsto \pm(\alpha\cdot x - t)^{s-1}_{+} : \|\alpha\|_1 = 1, \; |t| \leq 1 \}.
\end{align*}
We then set $ \mathcal{F}_{s} $ to be the linear span of $ \mathcal{H}_{s} $. With this notation, $ \mathcal{F}_{\text{ramp}} = \mathcal{F}_{2} $.

The condition 
$ \int_{\mathbb{R}^d}\|\omega\|^s_1|\widetilde{f}(\omega)|d\omega < +\infty $ ensures that $ f^{\star} $ (corrected by a $ (s-1) $-th degree ridge polynomial) belongs to $ L_{1,\mathcal{H}_{s}} $ and $ \|f^{\star}\|_{\mathcal{H}_{s}} \leq v_{f^{\star},s} $. Functions with moderate variation are particularly closely approximated. Nevertheless, even when $ \|f^{\star}\|_{\mathcal{H}} $ is infinite, we express the trade-offs in approximation accuracy for consistently estimating functions in the closure of the linear span of $ \mathcal{H} $.

In what follows, we assume that the internal parameters have $\ell_1$ norm at most $ {v_0} $. Likewise, we assume that $ x \in [-1, 1]^d $ so that $ |\theta \cdot x | \leq \|\theta\|_1 \leq v_0 $. This control on the size of the internal parameters will be featured prominently throughout. In the case of spline activation functions, we are content with the assumption $ {v_0} = 1 $. Note that if one restricts the size of the domain and internal parameters (say, to handle polynomials), the functions $ h $ are still bounded and Lipschitz but with possibly considerably worse constants.

Suppose data $ \{(X_i, Y_i)\}_{i=1}^n $ are independently drawn from the distribution of $ (X,Y) $. To produce predictions of the real-valued response $ Y $ from its input $ X $, the target regression function $ f^{\star}(x) = \mathbb{E}[Y|X=x] $ is to be estimated. The function $ f^{\star} $ is assumed to be bounded in magnitude by a positive constant $ B $. We assume the noise $ \varepsilon = Y - f^{\star}(X) $ has moments (conditioned on $ X $) that satisfy a Bernstein condition with parameter $ \eta > 0 $. That is, we assume
\begin{equation*}
\mathbb{E}(|\varepsilon|^k|X) \leq \frac{1}{2}k!\eta^{k-2}\mathbb{V}(\varepsilon|X), \qquad k = 3,4,\dots,
\end{equation*}
where $ \mathbb{V}(\varepsilon|X) \leq \sigma^2 $. This assumption is equivalent to requiring that $ \mathbb{E}(e^{|\varepsilon|/\nu}|X) $ is uniformly bounded in $ X $ for some $ \nu > 0 $. A stricter assumption is that $ \mathbb{E}(e^{|\varepsilon|^2/\nu}|X) $ is uniformly bounded in $ X $, which corresponds to an error distribution with sub-Gaussian tails. These two noise settings will give rise to different risk bounds, as we will see.


Because $ f^{\star} $ is bounded in magnitude by $ B $, it is useful to truncate an estimator $ \hat{f} $ at a level $ B_n $ at least $ B $. Depending on the nature of the noise $ \varepsilon $, we will see that $ B_n $ will need to be at least $ B $ plus a term of order $ \sqrt{\log n} $ or $ \log n $. We define the truncation operator $ T $ that acts on function $ f $ in $ \mathcal{F} $ by $ Tf = \min\{|f|,B_n\}\text{sgn}f $. Associated with the truncation operator is a tail quantity $
T_n = 2\sum_{i=1}^n(|Y_i|^2- B_n^2)\mathbb{I}\{ |Y_i| > B_n \} $ that appears in the following analysis and our risk bounds have a $ \mathbb{E}[T_n/n] $ term, but this will be seen to be negligible when compared to the main terms. The behavior of $ \mathbb{E}T_n $ is studied in \prettyref{lmm:truncate_bound}.

The empirical mean squared error of a function $ f $ as a candidate fit to the observed data is $ (1/n)\sum_{i=1}^n(Y_i-f(X_i))^2 $. Given the collection of functions $ \mathcal{F} $, a penalty $ \text{pen}_n(f) $, $ f \in \mathcal{F} $, and data, a penalized least squares estimator $ \hat{f} $ arises by optimizing or approximately optimizing
\begin{equation} \label{eq:least_squares}
(1/n)\sum_{i=1}^n(Y_i-f(X_i))^2 + \text{pen}_n(f)/n.
\end{equation}
Our method of risk analysis proceeds as follows. Given a collection $ \mathcal{F} $ of candidate functions, we show that there is a countable approximating set $ \widetilde{\mathcal{F}} $ of representations $ \widetilde{f} $, variable-distortion, variable-complexity cover of $ \mathcal{F} $, and a complexity function $ L_n(\widetilde{f}) $, with the property that for each $ f $ in $ \mathcal{F} $, there is an $ \widetilde{f} $ in $ \widetilde{\mathcal{F}} $ such that $ \text{pen}_n(f) $ is not less than a constant multiple of $ \gamma_f L_n(\widetilde{f}) + \Delta_n(f,\widetilde{f}) $, where $ \gamma_f $ is a constant (depending on $ B $, $ v_f $, $ \sigma^2 $, and $ \eta $) and $ \Delta_n(f,\widetilde{f}) $ is given as a suitable empirical measure of distortion (based on sums of squared errors). The variable-distortion, variable-complexity terminology has its origins in \cite{Barron2008-3, Barron2008-4, Cheang1998}.
The task is to determine penalties such that an estimator $ \hat{f} $ approximately achieving the minimum of $ \|Y-f\|^2_n + \text{pen}_n(f)/n $ satisfies
\begin{equation} \label{eq:risk}
\mathbb{E}\|T\hat{f}-f^{\star}\|^2 \leq c\inf_{f\in\mathcal{F}}\{ \|f-f^{\star}\|^2 + \mathbb{E}\text{pen}_n(f)/n \},
\end{equation}
for some universal $ c > 1 $. Valid penalties take different forms depending on the size of the effective dimension $ d $ relative to the sample size $ n $ and smoothness assumption of $ \phi $.

\begin{itemize}
\item When $ d $ is large compared to $ n $ and if $ \phi $ satisfies \prettyref{ass:ass1}, a valid penalty divided by sample size $ \text{pen}_n(f)/n $ is at least
\begin{equation} \label{eq:pen1}
16v_f\left(\frac{\gamma_n B_n^2{v^2_0}\log(d+1)}{n}\right)^{1/4} + 8\left(\frac{\gamma_n B_n^2{v^2_0}\log(d+1)}{n}\right)^{1/2} + \frac{T_n}{n}.
\end{equation}
\item When the noise $ \varepsilon $ is zero and $ d $ is large compared to $ n $ and if $ \phi $ satisfies \prettyref{ass:ass1}, a valid penalty divided by sample size $ \text{pen}_n(f)/n $ is at least
\begin{equation} \label{eq:pen2}
16v_f^{4/3}\left(\frac{\gamma_n{v^2_0}\log(d+1)}{n}\right)^{1/3}+ 4(v_f^{4/3}+1)\left(\frac{\gamma_n{v^2_0}\log(d+1)}{n}\right)^{2/3}.
\end{equation}
\item When $ d $ is large compared to $ n $ and if $ \phi $ satisfies \prettyref{ass:ass2}, a valid penalty divided by sample size $ \text{pen}_n(f)/n $ is at least of order
\begin{equation} \label{eq:pen3}
v_f^{4/3}\left(\frac{\gamma_n{v^2_0}\log(d+1)}{n}\right)^{1/3} + v_f\left(\frac{1}{n}\sum_{i=1}^n|Y_i|\right)\left(\frac{\gamma_n{v^2_0}\log(d+1)}{n}\right)^{1/3} +  \frac{T_n}{n}.
\end{equation}
\item When $ d $ is small compared to $ n $ and if $ \phi $ satisfies \prettyref{ass:ass1}, a valid penalty divided by sample size $ \text{pen}_n(f)/n $ is at least
\begin{align}
60v_f{v_0}\left(\frac{d\gamma_n\log(n/d+1)}{n}\right)^{1/2+1/2(d+3)}+\frac{1}{{v^2_0}}\left(\frac{d\gamma_n\log(n/d+1)}{n}\right)^{1/2+1/2(d+3)} \nonumber & \\ 
+ \left(\frac{d\gamma_n\log(n/d+1)}{n}\right)^{1/2+3/2(d+3)} + \frac{d\gamma_n\log(n/d+1)}{n} + \frac{T_n}{n}. \label{eq:pen4}
\end{align}
\end{itemize}

Here $ \gamma_n = (2\tau)^{-1}(1+\delta_1/2)(1+2/\delta_1)(B+B_n)^2 + 2(1+1/\delta_2)\sigma^2 + 2(B+B_n)\eta $ and $ \tau = (1+\delta_1)(1+\delta_2) $ for some $ \delta_1 > 0 $ and $ \delta_2 > 0 $.

Accordingly, if $ f^{\star} $ belongs to $ L_{1,\mathcal{H}} $, then $ \mathbb{E}\|T\hat{f}-f^{\star}\|^2 $ is not more than a constant multiple of the above penalties with $ v_f $ replaced by $ \|f^{\star}\|_{\mathcal{H}} $.

In the single-hidden layer case, we have the previously indicated quantification of the error of approximation $ \| f - f^{\star} \|^2 $. Nevertheless, the general result \prettyref{eq:risk} allows us to likewise say that the risk for multilayer networks will be at least as good as the deep network approximation capability will permit.
The quantity
\begin{equation*}
\inf_{f\in\mathcal{F}}\{ \|f-f^{\star}\|^2 + \mathbb{E}\text{pen}_n(f)/n \}.
\end{equation*}
is an index of resolvability of $ f^{\star} $ by functions $ \mathcal{F} $ with sample size $ n $. We shall take particular advantage of such risk bounds in the case that $ \text{pen}_n(f) $ does not depend on $ \underline{X} $. Our restriction of $ X $ to $ [-1,1]^d $ is one way to allow the construction of such penalties.

The following table expresses the heart of our results in the case of penalty based on the $ \ell_1 $ norm of the outer layer coefficients of one-hidden layer networks expressible through $ v_f $ (subject to constraints on the inner layer coefficients). These penalties also provide risk bounds for moderate and high-dimensional situations.

\begin{table}[H]
\centering
\caption{Main contributions to penalties for \\ \prettyref{thm:main_result} over continuum of candidate fits}
\label{tab:T1}
\begin{tabular}{c c c c c}
\hline
\hline
& Activation $ \phi $ & $ \text{pen}_n(f)/n \gtrsim $ & $ \lambda_n \gtrsim $ & \\ [0.5ex]
\hline
\hline
I & \prettyref{ass:ass1} & $ v_f\lambda_n $ & $ \left( \frac{\gamma^2_n\log(d+1)}{n} \right)^{1/4} $ & \\
\hline
II & \prettyref{ass:ass2} & $ (v_f)^{4/3}\lambda_n  $ & $ \left(\frac{\gamma^2_n\log(d+1)}{n} \right)^{1/3} $ & \\
\hline
III & \prettyref{ass:ass1} & $ v_f\lambda_n $ & $ \left( \frac{d\gamma_n\log(n/d+1)}{n} \right)^{1/2+1/(2(d+3))} $ & \\
\hline
\end{tabular}
\end{table}

\begin{table}[H]
\centering
\caption{Main contributions to penalties for \\ \prettyref{thm:main_result} over discretization of candidate fits}
\label{tab:T2}
\begin{tabular}{c c c c c}
\hline
\hline
& Activation $ \phi $ & $ \text{pen}_n(f)/n \gtrsim $ & $ \lambda_n \gtrsim $ & \\ [0.5ex]
\hline
\hline
A & \prettyref{ass:ass1} & $ (v_f)^{4/3}\lambda_n $ & $ \left( \frac{\gamma^2_n\log(d+1)}{n} \right)^{1/3} $ & \\
\hline
B & \prettyref{ass:ass2} & $ (v_f)^{6/5}\lambda_n  $ & $ \left(\frac{\gamma^2_n\log(d+1)}{n} \right)^{2/5} $ & \\
\hline
C & \prettyref{ass:ass1} & $ v_f\lambda_n $ & $ \left( \frac{d\gamma_n\log(n/d+1)}{n} \right)^{1/2+1/(2(d+3))} $ & \\
\hline
\end{tabular}
\end{table}

The results we wish to highlight are contained in the first two rows of \prettyref{tab:T1}. The penalties as stated are valid up to modest universal constants and negligible terms that do not depend on the candidate fit. The quantity $ \gamma_n $ is of order $ \log^2 n $ in the sub-exponential noise case, order $ \log n $ in the sub-Gaussian noise case and of constant order in the zero noise case. This $ \gamma_n $ (as defined in \prettyref{lmm:truncate_bound}) depends on the variance bound $ \sigma^2 $, Bernstein parameter $ \eta $, the upper bound $ B $ of $ \|f^{\star}\|_{\mathcal{H}} $, and the noise tail level $ B_n $ of the indicated order.

When $ f^{\star} $ belongs to $ L_{1,\mathcal{H}} $, a resulting valid risk bound is a constant multiple of $ \|f^{\star}\|_{\mathcal{H}}\lambda_n $ or $ \|f^{\star}\|^{4/3}_{\mathcal{H}}\lambda_n $, according to the indicated cases.
In this way the $ \lambda_n $ expression provides a rate of convergence. Thus the columns of \prettyref{tab:T1} provide valid risk bounds for these settings.

The classical risk bounds for mean squared error, involving $ d/n $ to some power, are only useful when the sample size is much larger than the dimension. Here, in contrast, in the first two lines of \prettyref{tab:T1}, we see the dependence on dimension is logarithmic, permitting much smaller sample sizes. These results are akin to those obtained in \cite{Wainwright2011} (where the role of the dimension there is the size of the dictionary) for high-dimensional linear regression. However, there is an important difference. Our dictionary of non-linear parameterized functions is infinite dimensional. For us, the role of $ d $ is the input dimension, not the size of the dictionary. The richness of $ L_{1,\mathcal{H}} $ is largely determined by the sizes of $ v_0 $ and $ v_f $ and $ L_{1,\mathcal{H}} $ more flexibly represents a larger class of functions.

The price we pay for the smaller dependence on input dimension is a deteriorated rate with exponent $ 1/4 $ in general and $ 1/3 $ under slightly stronger smoothness assumptions on $ \phi $, rather than the familiar exponents of $ 1/2 $.

The rate in the last row improves upon the familiar exponent of $ 1/2 $ to $ 1/2 + 1/(2(d+3)) $. Note that when $ d $ is large, this enhancement in the exponent is negligible. The rate in the first row is better than the third approximately for $ d > \sqrt{n} $, the second is better than the third row approximately for $ d > n^{1/3} $, and both of these first two rows have risk tending to zero as long as $ d < e^{o(n)} $. 

For functions in $ L_{1,\mathcal{H}_{\text{ramp}}} $, an upper bound of $ ((d/n)\log(n/d))^{1/2} $ for the squared error loss is obtained in \cite{Barron1994}. The $ L^2 $ squared error minimax rates for functions in $ L_{1,\mathcal{H}_{1}} = L_{1,\mathcal{H}_{\text{step}}} $ \cite{Barron1999}, was determined to be between 
\begin{equation*}
(1/n)^{1/2+1/(2(d+1))}(\log n)^{-(1+1/d)(1+2/d)(1+2/d)(2+1/d)}5
\end{equation*}
and
\begin{equation*}
(\log n/n)^{1/2+1/(2(2d+1))}.
\end{equation*}
 Using the truncated penalized $ \ell_1 $ least squares estimator \prettyref{eq:risk}, we obtain an improved rate of order $ ((d\gamma_n/n)\log(n/d))^{1/2+1/(2(d+3))} $, where $ \gamma_n $ is logarithmic in $ n $, using techniques that originate in \cite{Barron2008} and \cite{Huang2008}, with some corrections here. 

\section{How far from optimal?}
For positive $ v_0 $, let
\begin{equation}
\mathcal{D}_{v_0} = \mathcal{D}_{v_0, \phi} = \{ \phi(\theta \cdot x - t),\; x \in B: \|\theta\|_1 \leq v_0,\; t\in\mathbb{R} \}
\end{equation}
be the dictionary of all such inner layer ridge functions $ \phi(\theta \cdot x - t) $ with parameter restricted to the $ \ell_1 $ ball of size $ v_0 $ and variables $ x $ restricted to the cube $ [-1,1]^d $. The choice of the $ \ell_1 $ norm on the inner parameters is natural as it corresponds to $ \|\theta\|_B = \sup_{x\in B}|\theta\cdot x| $ for $ B = [-1,1]^d $. Let $ \mathcal{F}_{v_0, v_1} = \mathcal{F}_{v_0, v_1, \phi} = \ell_1(v_1, \mathcal{D}_{v_0}) $ be the closure of the set of all linear combinations of functions in $ \mathcal{D}_{v_0} $ with $ \ell_1 $ norm of outer coefficients not more than $ v_1 $. For any class of functions $ \mathcal{F} $ on $ [-1,1]^d $, the minimax risk is 
\begin{equation} \label{eq:minimax}
R_{n,d}(\mathcal{F}) = \inf_{\hat{f}} \sup_{f\in \mathcal{F}} \mathbb{E}\|f-\hat{f}\|^2,
\end{equation}

Consider the model $ Y = f(X) + \varepsilon $ for $ f \in \mathcal{F}_{v_0,v_1, \text{sine}} $, where $ \varepsilon \sim N(0, 1) $ and $ X \sim \text{Uniform}([-1,1]^d) $. It was determined in \cite{Klusowski2017}, that for $
\frac{d}{v_0}+1 > \left(c\frac{v^2_1n}{v_0\log(1+d/v_0)}\right)^{1/v_0} $, roughly corresponding to $ d \gg n $,
\begin{equation} \label{eq:lowerhighdim}
R_{n,d}(\mathcal{F}_{v_0,v_1, \text{sine}}) \geq C\left(\frac{v_0v^2_1\log(1+d/v_0)}{n}\right)^{1/2},
\end{equation}
and for $ \frac{v_0}{d}+1 > \left(c\frac{v^2_1n}{d\log(1+v_0/d)}\right)^{1/d} $,
\begin{equation} \label{eq:lowerlowdim}
R_{n,d}(\mathcal{F}_{v_0,v_1, \text{sine}}) \geq C\left(\frac{dv^2_1\log(1+v_0/d)}{n}\right)^{1/2}.
\end{equation}

These lower bounds are similar in form to the risk upper bounds that are implied from the penalties in \prettyref{tab:T2}. These quantities have the attractive feature that the rate (the power of $ 1/n $) remains at least as good as $ 1/2 $ or $ 2/5 $ even as the dimension grows. However, rates determined by \prettyref{eq:lowerlowdim} and the last line in \prettyref{tab:T2} are only useful provided $ d/n $ is small. In high dimensional settings, the available sample size might not be large enough to ensure this condition. 

These results are all based on obtaining covering numbers for the library $ \{ x\mapsto \phi(\theta\cdot x) : \|\theta\|_1 \leq {v_0} \} $. If $ \phi $ satisfies a Lipschitz condition, these numbers are equivalent to $\ell_1$ covering numbers of the internal parameters or of the Euclidean inner product of the data and the internal parameters. The factor of $ d $ multiplying the reciprocal of the sample size is produced from the order $ d\log({v_0}/\epsilon) $ log cardinality of the standard covering of the library $ \{\theta:\|\theta\|_1 \leq {v_0} \} $. What enables us to circumvent this polynomial dependence on $ d $ is to use an alternative cover of $ \{x\mapsto x\cdot \theta: \|\theta\|_1 \leq {v_0} \} $ that has log cardinality of order $ ({v_0}/\epsilon)^2\log(d+1) $. Misclassification errors for neural networks with bounded internal parameters have been analyzed in \cite{Bartlett1998, Bartlett2009, Lee1996} (via Vapnik-Chervonenkis dimension and its implications for covering numbers). Unlike the setup considered here, past work \cite{Barron1994, Barron1999-2, Barron1999, Barron2008, Barron2008-2, Barron2008-4, Lee1996, Bunea2007, Martin1999, Zhang2009, Nemirovski2000} has not investigated the role of such restricted parameterized classes in the determination of suitable penalized least squares criterion for non-parametric function estimation. After submission of the original form of this work, our results have been put to use in \cite{Schmidt-Hieber2017} to give risk statements about multi-layer (deep) networks activated by ramp functions.

\section{Computational aspects}

From a computational point of view, the empirical risk minimization problem \prettyref{eq:least_squares} is highly non-convex, and it is unclear why existing algorithms like gradient descent or back propagation are empirically successful at learning the representation \prettyref{eq:general}.  There are relatively few rigorous results that guarantee learning for regression models with latent variables, while keeping both the sampling and computational complexities polynomial in $ n $ and $ d $. Here we catalogue some papers that make progress toward developing a provably good, computationally feasible estimation procedure. Most of them deal with parameter recovery and assume that $ f^{\star} $ has exactly the form \prettyref{eq:general}. Using a theory of tensor decompositions from \cite{Anandkumar2014}, \cite{Janzamin2015} apply the method of moments via tensor factorization techniques to learn mixtures of sigmoids, but they require a special non-degeneracy condition on the activation function. It is assumed that the input distribution $ P $ is known apriori. In \cite{Bartlett2017}, the authors use tensor initialization and resampling to learn the parameters in a representation of the form \prettyref{eq:general} with smooth $ \phi $ that has sample complexity $ O(d) $ and computation complexity $ O(dn) $.

In \cite{Montanari2017}, the authors estimate the gradient of the regression function (where $ X $ is Gaussian and $ \phi $ is the logistic sigmoid) at a set of random points, and then cluster the estimated gradients. They prove that the estimated gradients concentrate around the internal parameter vectors. However, unless the weights of the outer layer are positive and sum to $ 1 $, the complexity is exponential in $ d $. In \cite{Valiant2014}, it was shown that for a randomly initialized neural network with sufficiently many hidden units, the generic gradient descent algorithm learns any low degree polynomial. Learning non-linear networks through multiple rounds of random initialization followed by arbitrary optimization steps was proposed in \cite{Zhang2016}. In \cite{Zhang2015}, an efficiently learned kernel based estimator was shown to perform just as well as a class of deep neural networks. However, its ability to well-approximate general conditional mean regression functions is unclear.

The next section discusses an iterative procedure that reduces the complexity of finding the penalized least squares estimator \prettyref{eq:least_squares}.

\section{Greedy algorithm}

The main difficulty with constructing an estimator that satisfies \prettyref{eq:risk} is that it involves a $ dm $-dimensional optimization.
Here, we outline a greedy approach that reduces the problem to performing $ m $ $ d $-dimensional optimizations. This construction is based on the $\ell_1$-penalized greedy pursuit (LPGP) in \cite{Barron2008}, with the modification that the penalty can be a convex function of the candidate function complexity. Greedy strategies for approximating functions in the closure of the linear span of a subset of a Hilbert space has its origins in \cite{Lee1992} and many of its statistical implications were studied in \cite{Barron2008-2} and \cite{Barron2008}.

Let $ f^{\star} $ be a function, not necessarily in $ \mathcal{F} $. Initialize $ f_0 = 0 $. For $ m = 1,2,\dots, $ iteratively, given the terms of $ f_{m-1} $ as $ h_1,\dots,h_{m-1} $ and the coefficients of it as $ \beta_{1,m-1},\dots,\beta_{m-1,m-1} $, we proceed as follows. Let $ f_m(x) = \sum_{j=1}^{m}\beta_{j,m}h_{j}(x) = \sum_{j=1}^{m}\beta_{j,m}\phi(\theta_{h_{j}}\cdot x) $, with the term $ h_{m} $ in $ \mathcal{H} $ chosen to come within a constant factor $ c \geq 1 $ of the maximum inner product with the residual $ f^{\star} - f_{m-1} $; that is
\begin{equation*}
 \langle h_{m}, f^{\star} - f_{m-1} \rangle \geq \frac{1}{c}\sup_{h\in\mathcal{H}}\langle h, f^{\star} - f_{m-1} \rangle.
\end{equation*}
Define $ f_m(x) = (1-\alpha_m)f_{m-1}(x) + \beta_{m,m}h_{m}(x) $. Associated with this representation of $ f_m $ is the $ \ell_1 $ norm of its coefficients $ v_m = \sum_{j=1}^m|\beta_{j,m}| = (1-\alpha_m)v_{m-1} + \beta_{m,m} $.
The coefficients $ \alpha_m $ and $ \beta_{m,m} $ are chosen to minimize $ \|f^{\star}-(1-\alpha_m)f_{m-1} - \beta_{m,m}h_{m} \|^2 + \omega((1-\alpha_m)v_{m-1} + \beta_{m,m}) $.

In the empirical setting, with $ R_i = Y_i - f_{m-1}(X_i) $, the high-dimensional optimization task is to find $ \theta_m $ such that
\begin{equation*}
\frac{1}{n}\sum_{i=1}^n R_i\phi(\theta_m\cdot X_i) \geq \frac{1}{c}\sup_{\theta}\frac{1}{n}\sum_{i=1}^n R_i\phi(\theta\cdot X_i)
\end{equation*}
The fact that one does not need to find the exact maximizer of the above empirical inner product, but only come to within a constant multiple of it, has important consequences. For example, in adaptive annealing, one begins by sampling from an initial distribution $ p_0 $ and then iteratively samples from a distribution proportional to $ e^{t(\frac{1}{n}\sum_{i=1}^n R_i\phi(\theta\cdot X_i))}p_0(\theta) $, evolving according to $ \theta_{t+h} = \theta_t - hG_t(\theta_t) $, where $ G_t(\theta) $ satisfies $ \nabla^T [G_t(\theta)p_t(\theta)]  = \partial_t p_t(\theta)  $. The mean of $ p_t $ is at least $ \frac{1}{c}\sup_{\|\theta\|_1 \leq \Lambda}\frac{1}{n}\sum_{i=1}^n R_i\phi(\theta\cdot X_i) $ for sufficiently large $ t $.

\begin{theorem} \label{thm:greedy}
Suppose $ w:\mathbb{R}\rightarrow \mathbb{R} $ is a real-valued non-negative convex function. If $ f_m $ is chosen according to the greedy scheme described previously, then
\begin{equation} \label{eq:greedy1}
\|f^{\star}-f_m\|^2 + w(v_m) \leq \inf_{f\in\mathcal{F}} \left\{ \|f^{\star} - f \|^2 + w(cv_f) + \frac{4b_f}{m}\right\},
\end{equation}
where $ b_f = c^2v^2_f + 2v_f\|f^{\star}\|(c+1) - \|f\|^2 $.
Furthermore, for all $ \delta > 0 $,
\begin{align} \label{eq:greedy2}
& \nonumber \|f^{\star}-f_m\|^2 + w(v_m) \\ & \leq \inf_{f\in\mathcal{F}}\inf_{\delta > 0}\left\{ (1+\delta)\|f^{\star} - f \|^2 + w(cv_f) + \frac{4(1+\delta)\delta^{-1}(c+1)^2v^2_f}{m}\right\},
\end{align}
and hence with $ \delta = \frac{2(c+1)v_f}{\|f^{\star}-f\|\sqrt{m}} $,
\begin{equation*}
\|f^{\star}-f_m\|^2 + w(v_m) \leq \inf_{f\in\mathcal{F}}\left\{ \left(\|f^{\star} - f \| + \frac{2(c+1)v_f}{\sqrt{m}}\right)^2 + w(cv_f)\right\}.
\end{equation*}
\end{theorem}

\begin{proof}
Fix any $ f $ in the linear span $ \mathcal{F} $, with the form $ \sum_{h\in\mathcal{H}}\beta_h h $, with non-negative $ \beta_h $ and set
\begin{equation*}
e_m = \|f^{\star} - f_m \|^2 - \|f^{\star}-f\|^2 + w(v_m).
\end{equation*}
From the definition of $ \alpha_m $ and $ \beta_{m,m} $ as minimizers of $ e_m $ for each $ h_m $, and the convexity of $ w $,
\begin{align*}
e_m 
& = \|f^{\star} - (1-\alpha_m)f_{m-1} - \beta_{m,m}h_{m} \|^2 - \|f^{\star}-f\|^2 + \\ & \qquad w((1-\alpha_m)v_{m-1}+\beta_{m,m}) \\
& \leq \|f^{\star} - (1-\alpha_m)f_{m-1} - \alpha_mcv_fh_m \|^2 - \|f^{\star}-f\|^2 + \\ & \qquad w((1-\alpha_m)v_{m-1}+\alpha_m cv_f) \\
& \leq \|f^{\star} - (1-\alpha_m)f_{m-1} - \alpha_mcv_fh_m \|^2 - \|f^{\star}-f\|^2 + \\ & \qquad (1-\alpha_m)w(v_{m-1})+\alpha_mw(cv_f).
\end{align*}
Now $ \|f^{\star} - (1-\alpha_m)f_{m-1} - \alpha_mcv_fh_m \|^2 $ is equal to $ \| (1-\alpha_m)(f^{\star} - f_{m-1}) + \alpha_m(f^{\star}- ch_mv_f) \|^2 $. Expanding this quantity leads to
\begin{align*}
\|f^{\star} - (1-\alpha_m)f_{m-1} - \alpha_mcv_fh_m \|^2 
& = (1-\alpha_m)^2\|f^{\star} - f_{m-1}\|^2 \\
& - 2\alpha_m(1-\alpha_m)\langle f^{\star} - f_{m-1}, ch_mv_f - f^{\star} \rangle \\ & +  \alpha^2_m\|f^{\star} - ch_mv_f\|^2.
\end{align*}
Next we add $ (1-\alpha_m)w(v_{m-1})+\alpha_m w(cv_f) - \|f^{\star}-f\|^2 $ to this expression to obtain
\begin{align}
e_m
& \leq (1-\alpha_m)e_{m-1} + \alpha^2_m[\|f^{\star}-ch_mv_f\|^2-\|f^{\star}-f\|^2] + \alpha_mw(cv_f) \nonumber \\
& \qquad - 2\alpha_m(1-\alpha_m)\langle f^{\star} - f_{m-1}, ch_mv_f - f \rangle \nonumber \\ & \qquad\qquad + \alpha_m(1-\alpha_m)[2\langle f^{\star}-f_{m-1}, f^{\star} - f\rangle - \|f^{\star}-f_{m-1}\|^2 - \|f^{\star}-f\|^2] \label{eq:norm}.
\end{align}
The expression in brackets in \prettyref{eq:norm} is equal to $ -\|f-f_{m-1}\|^2 $ and hence the entire quantity is further upper bounded by
\begin{align*}
e_m 
& \leq (1-\alpha_m)e_{m-1} + \alpha^2_m[\|f^{\star}-ch_mv_f\|^2-\|f^{\star}-f\|^2] + \alpha_mw(cv_f) \\
& \qquad - 2\alpha_m(1-\alpha_m)\langle f^{\star} - f_{m-1}, ch_mv_f - f \rangle.
\end{align*}
Consider a random variable that equals $ h $ with probability $ \beta_h/v_f $ having mean $ f $. Since a maximum is at least an average, the choice of $ h_m $ implies that $ \langle f^{\star} - f_{m-1}, ch_mv_f \rangle $ is at least $ \langle f^{\star} - f_{m-1}, f \rangle $. This shows that $ e_m $ is no less than $ (1-\alpha_m)e_{m-1} + \alpha^2_m[\|f^{\star}-ch_mv_f\|^2-\|f^{\star}-f\|^2] + \alpha_m w(cv_f) $. Expanding the squares in $ \|f^{\star}-ch_mv_f\|^2-\|f^{\star}-f\|^2 $ and using the Cauchy-Schwarz inequality yields the bound $ \|ch_mv_f\|^2 + 2\|f^{\star}\|(\|f-ch_mv_f\|) - \|f\|^2 $. Since $ \|h_m\| \leq \|h_m\|_{\infty} \leq 1 $ and $ \|f\| \leq v_f $, we find that $ \|f^{\star}-ch_mv_f\|^2-\|f^{\star}-f\|^2 $ is at most $ b_f = c^2v^2_f + 2v_f\|f^{\star}\|(c+1) - \|f\|^2 $.
Hence we have shown that
\begin{equation*}
e_1 \leq b_f + w(cv_f)
\end{equation*}
and
\begin{equation} \label{eq:greedy}
e_m \leq (1-\alpha_m)e_{m-1} + \alpha^2_mb_f + \alpha_mw(cv_f).
\end{equation}
Because $ \alpha $ is a minimizer of $ e_m $, it can replace it by any value in $ [0, 1] $ and the bound \prettyref{eq:greedy} holds verbatim. In particular, we can choose $ \alpha_m = 2/(m+1) $, $ m \geq 2 $ and use an inductive argument to establish \prettyref{eq:greedy1}.
The second statement \prettyref{eq:greedy2} follows from similar arguments upon consideration of 
\begin{equation*}
e_m = \|f^{\star} - f_m \|^2 - (1+\delta)\|f^{\star}-f\|^2 + w(v_m),
\end{equation*}
together with the inequality $ a^2-(1+\delta)b^2 \leq (1+\delta)\delta^{-1}(a-b)^2 $.
\end{proof}

\section{Risk bounds}

\subsection{Penalized estimators over a discretization of the parameter space}

In the case that $ \mathcal{F}_{\epsilon} $ is an $ L^2(P) $ $ \epsilon $-net of the parameter space and the noise is bounded, it follows from \prettyref{thm:fundamental} and \prettyref{eq:risk_bound} that if $ \text{pen}_n(f) = \gamma _f L_n(f) = \gamma _n \log \text{card}(\mathcal{F}_{\epsilon}) $, then
\begin{align*}
\mathbb{E}\|\hat{f}-f^{\star}\|^2 & \leq (\tau+1)\inf_{f\in\mathcal{F}_{\epsilon}} \{ \|f-f^{\star}\|^2 + \gamma _f \text{card}(\mathcal{F}_{\epsilon})/n \} \\
& \leq (\tau+1)\inf_{\epsilon > 0} \{ \epsilon^2 + \gamma _{f^{\star}} \log \text{card}(\mathcal{F}_{\epsilon})/n \}.
\end{align*}
By Theorem in \cite{Mendelson2002}, there exists a universal constant $ C_{{v_0}} > 0 $ such that $ \log \text{card}(\mathcal{F}_{\epsilon}) \leq C_{{v_0}}d\epsilon^{-\frac{2d}{d+2}} $. Hence,
\begin{align*}
\mathbb{E}\|\hat{f}-f^{\star}\|^2 & \leq (\tau+1)\inf_{f\in\mathcal{F}_{\epsilon}} \{ \|f-f^{\star}\|^2 + \gamma _f \text{card}(\mathcal{F}_{\epsilon})/n \} \\
& \leq (\tau+1)\inf_{\epsilon > 0} \{ \epsilon^2 + \gamma _{f^{\star}} C_{{v_0}} d\epsilon^{-\frac{2d}{d+2}}/n \} \\
& \leq (\tau+1)\left( \frac{C_{{v_0}} \gamma _{f^{\star}}d}{n} \right)^{\frac{d+2}{2d+2}}.
\end{align*}
This result is similar to \cite{Xiaohong1999}, which also improved on the more familiar rate of $ (\frac{d}{n})^{1/2} $ are obtained.

On the other hand, if $ h = \phi(x\cdot \theta_h) $, $ \phi $ satisfies \prettyref{ass:ass2}, and $ \|\theta_h \|_1 \leq {v_0} $, we can use an alternative argument via \prettyref{lmm:approx2} to produce $ \log \text{card}(\mathcal{F}_{\epsilon}) \leq C\epsilon^{-3}v^3v^2_0\log(d+1) $ for functions with variation at most $ v $. Hence,
\begin{align*}
\mathbb{E}\|\hat{f}-f^{\star}\|^2 & \leq (\tau+1)\inf_{f\in\mathcal{F}_{\epsilon}} \{ \|f-f^{\star}\|^2 + \gamma _n  \text{card}(\mathcal{F}_{\epsilon})/n \} \\
& \leq (\tau+1)\inf_{\epsilon > 0} \{ \epsilon^2 +  \epsilon^{-3}v^3_f v_0^2\log(d+1)/n \} \\
& \leq (\tau+1)\left( \frac{C \gamma_n v^3_{f^{\star}} v^2_0\log(d+1)}{n} \right)^{2/5}.
\end{align*}
Compare this result with the minimax risk lower bound \prettyref{eq:lowerhighdim} of order $ (\frac{\log(d+1)}{n})^{1/2} $. This conclusion for the discretized parameter space holds for any bounded Lipschitz activation function. In the extension to optimize over the continuum in the next section, we obtain the $ 1/3 $ power rate only under the stronger \prettyref{ass:ass2} and a $ 1/4 $ rate for the general bounded Lipschitz case. This attainment of a rate for function estimation including single-hidden layer neural nets which is a power of $ (\log d) /n $ rather than $ d/n $ is the heart of the novel contribution of this paper.

\subsection{Penalized estimators over the entire parameter space}

Here we state our main theorem.
\begin{theorem} \label{thm:main_result}
Let $ f^{\star} $ be a real-valued function on $ [-1,1]^d $ with finite variation $ v_{f^{\star}} $ with respect to the library $ \mathcal{H} = \{  h(x) = \phi(\theta \cdot x) : \|\theta\|_1 \leq v_0 \} $.
If $ \hat{f} $ is chosen to satisfy
\begin{equation*}
\frac{1}{n}\sum_{i=1}^n(Y_i-\hat{f}(X_i))^2 + \text{pen}_n(\hat{f})/n \leq \inf_{f\in\mathcal{F}} \left\{ \frac{1}{n}\sum_{i=1}^n(Y_i-f(X_i))^2 + \text{pen}_n(f)/n \right\},
\end{equation*}
then for the truncated estimator $ T\hat{f} $ and for $ \text{pen}_n(f) $ depending on $ v_f $ as specified below, the risk has the resolvability bound
\begin{equation*}
\mathbb{E}\|T\hat{f}-f^{\star}\|^2 \leq (\tau+1)\inf_{f\in\mathcal{F}}\{ \|f-f^{\star}\|^2 + \mathbb{E}\text{pen}_n(f)/n \},
\end{equation*}
with penalfties as described in \prettyref{eq:pen1}, \prettyref{eq:pen2}, \prettyref{eq:pen3}, and \prettyref{eq:pen4}.
If $ \hat{f}_m $ is the LPGP estimator from the previous section, then by \prettyref{thm:greedy},
\begin{equation*}
\frac{1}{n}\sum_{i=1}^n(Y_i-\hat{f}_m(X_i))^2 + w(v_{\hat{f}_m}) \leq \inf_{f\in\mathcal{F}} \left\{ \frac{1}{n}\sum_{i=1}^n(Y_i-f(X_i))^2 + w(cv_f) + \frac{4b_f}{m} \right\},
\end{equation*}
where $ b_f $ is the empirical version of the same quantity in \prettyref{thm:greedy} and hence
the risk has the resolvability bound
\begin{equation*}
\mathbb{E}\|T\hat{f}-f^{\star}\|^2 \leq (\tau+1)\inf_{f\in\mathcal{F}}\{ \|f-f^{\star}\|^2 + \mathbb{E}\text{pen}_n(cf)/n + 4\mathbb{E}b_f/m \},
\end{equation*}
for a penalty, convex in $ v_f $, $ \text{pen}_n(f) = nw(v_f) $ as before.
If $ m $ is chosen to be of order between $ \sqrt{n} $ and $ n $ so as to make the computational effects negligible, the previously described $ L^2(P) $ rates for estimating $ f^{\star} $ in $ L_{1,\mathcal{H}} $ via the truncated estimator $ T\hat{f}_m $ are attainable under the appropriate penalties.
\end{theorem}


One can also extend these results to include penalties that depend on the number of terms $ m $ in an $ m $-term greedy approximation $ \hat{f}_m $ to $ f^{\star} $. We take $ \hat{f}_m $ to be an $ m $ term fit from an LPGP algorithm and choose $ \hat{m} $ among all $ m \in \mathcal{M} $ (i.e. $ \mathcal{M} = \{1,\dots, n\} $) to minimize
\begin{equation*}
\frac{1}{n}\sum_{i=1}^n(Y_i - \hat{f}_m(X_i))^2 + \text{pen}_n(\hat{f}_m,m)/n.
\end{equation*}
This approach enables the use of a data-based stopping criterion for the greedy algorithm. For more details on these adaptive methods, we refer the reader to \cite{Barron2008}. The resolvability risk bound allows also for interpolation rates between $ L_2 $ and $ L_{1,\mathcal{H}} $ refining the results of \cite{Barron2008-2} and in accordance with the best balance between error of approximation and penalty.

The target $ f^{\star} $ is not necessarily in $ \mathcal{F} $. To each $ f $ in $ \mathcal{F} $, there corresponds a function $ \rho $, which assigns to $ (X,Y) $ the relative loss
\begin{align*}
\rho(X,Y) = \rho_f(X,Y) = (Y-f(X))^2 - (Y-f^{\star}(X))^2.
\end{align*}
Let $ (\underline{X}^{\prime}, \underline{Y}^{\prime}) $ be an independent copy of the training data $ (\underline{X}, \underline{Y}) $ used for testing the efficacy of a fit $ \hat{f} $ based on $ (\underline{X}, \underline{Y}) $. The relative empirical loss with respect to the training data is denoted by $ P_n(f||f^{\star}) = \frac{1}{n}\sum_{i=1}^n\rho(X_i,Y_i) $ and that with respect to the independent copy is $ P^{\prime}_n(f||f^{\star}) = \frac{1}{n}\sum_{i=1}^n\rho(X^{\prime}_i,Y^{\prime}_i) $. We define the empirical squared error on the training and test data by $ D_n(f,\widetilde{f}) = \frac{1}{n}\sum_{i=1}^n(f(X_i) - \widetilde{f}(X_i))^2 $ and $ D^{\prime}_n(f,\widetilde{f}) = \frac{1}{n}\sum_{i=1}^n(f(X^{\prime}_i) - \widetilde{f}(X^{\prime}_i))^2 $ for all $ f,\widetilde{f} $ in $ \mathcal{F} $. Using the relationship $ Y = f^{\star}(X) + \varepsilon $, we note that $ \rho(X,Y) $ can also be written as $ (f(X)-f^{\star}(X))^2 - 2\varepsilon(f(X)-f^{\star}(X)) = g^2(X) - 2\varepsilon g(X) $, where $ g(x) = f(x) - f^{\star}(x) $. Hence we have the relationship $ P_n(f||f^{\star}) = D_n(f,f^{\star})-\frac{2}{n}\sum_{i=1}^n\varepsilon_ig(X_i) $.

The relative empirical loss $ P^{\prime}_n(\hat{f}||f^{\star}) $ is an unbiased estimate of the risk $ \mathbb{E}\|\hat{f}-f^{\star}\|^2 $. Since $ \varepsilon^{\prime}_i $ has mean zero conditioned on $ X^{\prime}_i $, the mean of $ P^{\prime}_n(\hat{f}||f^{\star}) $ with respect to $ (\underline{X}^{\prime}, \underline{Y}^{\prime}) $ is $ \|\hat{f}-f^{\star}\|^2 $. This quantity captures how well the fit $ \hat{f} $ based on the training data generalizes to a new set of observations. The goal is to control the empirical discrepancy $ P^{\prime}_n(f||f^{\star}) - \tau P_n(f||f^{\star}) $ between the loss on the future data and the loss on the training data for a constant $ \tau > 1 $. Toward this end, we seek a positive quantity $ \text{pen}_n(f) $ to satisfy
\begin{equation*}
\mathbb{E}\sup_{f\in\mathcal{F}}\left\{ P_n^{\prime}(f||f^{\star}) - \tau P_n(f||f^{\star}) - \tau\text{pen}_n(f)/n \right\} \leq 0,
\end{equation*}
Once such an inequality holds, the data-based choice $ \hat{f} $ in $ \mathcal{F} $ yields
\begin{equation*}
\mathbb{E}P_n^{\prime}(\hat{f}||f^{\star}) \leq \tau\mathbb{E}[P_n(\hat{f}||f^{\star}) + \text{pen}_n(f)/n].
\end{equation*}
If $ \hat{f} $ satisfies
\begin{equation} \label{eq:estimator}
\frac{1}{n}\sum_{i=1}^n(Y_i-\hat{f}(X_i))^2 + \frac{\text{pen}_n(\hat{f})}{n} \leq \inf_{f\in\mathcal{F}} \left\{ \frac{1}{n}\sum_{i=1}^n(Y_i-f(X_i))^2 + \frac{\text{pen}_n(f)}{n} + A_f \right\},
\end{equation}
for some positive quantity $ A_f $ that decays to zero as the sample size grows,
we see that
\begin{equation*}
\mathbb{E}P_n^{\prime}(\hat{f}||f^{\star}) \leq \tau\inf_{f\in\mathcal{F}} \mathbb{E}[P_n(f||f^{\star}) + \text{pen}_n(f)/n + A_f].
\end{equation*}
Using $ \mathbb{E}P_n^{\prime}(\hat{f}||f^{\star}) = \mathbb{E}\|\hat{f}-f^{\star}\|^2 $ and $ \mathbb{E}P_n(f||f^{\star}) = \|f-f^{\star}\|^2 $, the above expression is seen to be
\begin{equation} \label{eq:risk_bound}
\mathbb{E}\|\hat{f}-f^{\star}\|^2 \leq \tau\inf_{f\in\mathcal{F}} \{ \|f-f^{\star}\|^2 + \mathbb{E}\text{pen}_n(f)/n + \mathbb{E}A_f \}.
\end{equation}
For the purposes of proving results in the case when $ \mathcal{F} $ is uncountable, it is useful to consider complexities $ L_n(\widetilde{f}) $ for $ \widetilde{f} $ in a countable subset $ \widetilde{\mathcal{F}} $ of $ \mathcal{F} $ satisfying $ \sum_{\widetilde{f}\in\widetilde{\mathcal{F}}}e^{-\gamma_n L_n(\widetilde{f})} \leq 1 $ for some $ \gamma_n > 0 $ and such that
\begin{align*}
\sup_{f\in\mathcal{F}}\left\{ P_n^{\prime}(f||f^{\star}) - \tau P_n(f||f^{\star}) - \tau\text{pen}_n(f)/n \right\}  & \\
\leq \sup_{\widetilde{f}\in\widetilde{\mathcal{F}}}\left\{ P_n^{\prime}(\widetilde{f}||f^{\star}) - \tau P_n(\widetilde{f}||f^{\star}) - \tau\gamma_n L_n(\widetilde{f})/n \right\},
\stepcounter{equation}\tag{\theequation}\label{eq:countable}
\end{align*}
with
\begin{align*}
\mathbb{E}\sup_{\widetilde{f}\in\widetilde{\mathcal{F}}}\left\{ P_n^{\prime}(\widetilde{f}||f^{\star}) - \tau P_n(\widetilde{f}||f^{\star}) - \tau\gamma_n L_n(\widetilde{f})/n \right\} \leq 0.
\end{align*}
The condition in \prettyref{eq:countable} is equivalent to requiring that
\begin{equation*}
\sup_{f\in\mathcal{F}}\inf_{\widetilde{f}\in\widetilde{\mathcal{F}}}\{ \Delta_n(f,\widetilde{f}) + \gamma_n L_n(\widetilde{f}) - \text{pen}_n(f)\} \leq 0,
\end{equation*}
where
\begin{align*} 
\Delta_n(f,\widetilde{f}) = n[P_n(\widetilde{f}||f^{\star})-P_n(f||f^{\star})] - (n/\tau)[P^{\prime}_n(\widetilde{f}||f^{\star}) - P^{\prime}_n(f||f^{\star})].
\end{align*}
If we truncate the penalized least squares estimator $ \hat{f} $ at a certain level $ B_n $, for $ \mathbb{E}\|T\hat{f}-f^{\star}\|^2 $ to maintain the resolvability bound $ \tau\inf_{f\in\mathcal{F}} \{ \|f-f^{\star}\|^2 + \mathbb{E}\text{pen}_n(f)/n + \mathbb{E}A_f \} $, we require that
\begin{equation*}
\sup_{f\in\mathcal{F}}\inf_{\widetilde{f}\in\widetilde{\mathcal{F}}}\{ \Delta_n(f,\widetilde{f}) + \gamma_n L_n(\widetilde{f}) - \text{pen}_n(f)\} \leq 0,
\end{equation*}
where
\begin{align*}
\Delta_n(f,\widetilde{f}) = 
n[P_n(T\widetilde{f}||f^{\star})-P_n(f||f^{\star})] - (n/\tau)[P^{\prime}_n(T\widetilde{f}||f^{\star}) - P^{\prime}_n(Tf||f^{\star})].
\end{align*}
Rather than working with the relative empirical loss $ P^{\prime}_n(Tf||f^{\star}) $, we prefer to work with $ D^{\prime}_n(Tf,f^{\star}) $. These two quantities are related to each other, provided $ \frac{1}{n}\sum_{i=1}^n\varepsilon_ig(X^{\prime}_i) $ is small and they are exactly equal in the no noise case. Hence we would like to determine penalties that ensure
\begin{equation*}
\mathbb{E}\sup_{f\in\mathcal{F}}\left\{ D_n^{\prime}(Tf,f^{\star}) - \tau P_n(f||f^{\star}) - \tau\text{pen}_n(f)/n \right\} \leq 0.
\end{equation*}
Suppose we require that
\begin{equation*}
\mathbb{E}\sup_{f\in\mathcal{F}}\{ \tau_1^{-1}D_n^{\prime}(Tf,f^{\star}) - \tau P_n(f||f^{\star}) - \tau\text{pen}_n(f)/n \} \leq 0,
\end{equation*}
for some $ \tau_1 \geq 1 $. This further inflates the resulting risk bound by $ \tau_1 $ so that the factor $ \tau $ is replaced with $ \tau\tau_1 $ in \prettyref{eq:risk_bound}. However, it enables us to create countable covers $ \widetilde{\mathcal{F}} $ with smaller errors in approximating functions from $ \mathcal{F} $. To see this, suppose the countable cover $ \widetilde{\mathcal{F}} $ satisfies
\begin{align*}
\sup_{f\in\mathcal{F}}\left\{ \tau_1^{-1}D_n^{\prime}(Tf,f^{\star}) - \tau P_n(f||f^{\star}) - \tau\text{pen}_n(f)/n \right\}  & \\ 
\leq \sup_{\widetilde{f}\in\widetilde{\mathcal{F}}}\left\{ D_n^{\prime}(T\widetilde{f},f^{\star}) - \tau P_n(T\widetilde{f}||f^{\star}) - \tau\gamma_n L_n(\widetilde{f})/n \right\},
\end{align*}
or equivalently that
\begin{equation*}
\sup_{f\in\mathcal{F}}\inf_{\widetilde{f}\in\widetilde{\mathcal{F}}}\left\{ \Delta_n(f,\widetilde{f}) + \gamma_n L_n(\widetilde{f}) - \text{pen}_n(f)\right\} \leq 0,
\end{equation*}
where
\begin{align*} 
\Delta_n(f,\widetilde{f})
& = n[P_n(T\widetilde{f}||f^{\star})- P_n(f||f^{\star})] + \\ & \qquad
n\tau^{-1}[\tau_1^{-1}D^{\prime}_n(Tf,f^{\star}) - D^{\prime}_n(T\widetilde{f},f^{\star})].
\end{align*}
We set $ \tau_1 = 1/\tau+1 $. Using the inequality, $ \tau^{-1}a^2- b^2 \leq \frac{1}{\tau-1}(b-a)^2 $ that can be derived from $ (a/\sqrt{\tau}-b\sqrt{\tau})^2 \geq 0 $, we can upper bound the difference $ \tau_1^{-1}D^{\prime}_n(Tf,f^{\star}) - D^{\prime}_n(T\widetilde{f},f^{\star}) $ by
\begin{align*} 
(\tau_1-1)^{-1}D^{\prime}_n(Tf,T\widetilde{f}).
\end{align*}
This quantity does not involve $ f^{\star} $, which is desirable for the proceeding analysis.
Hence $ \Delta_n(f,\widetilde{f}) $ is not greater than
\begin{align*} 
n[P_n(T\widetilde{f}||f^{\star})-P_n(f||f^{\star}) + D^{\prime}_n(Tf,T\widetilde{f})].
\end{align*}
and thus we seek a penalty $ \text{pen}_n(f) $ that is at least
\begin{align*} 
\inf_{\widetilde{f}\in\widetilde{\mathcal{F}}} \{ \gamma_n L_n(\widetilde{f}) + n[P_n(T\widetilde{f}||f^{\star})-P_n(f||f^{\star}) + D^{\prime}_n(Tf,T\widetilde{f})] \}.
\end{align*}
An estimator $ \hat{f} $ satisfying \prettyref{eq:estimator} with penality $ \text{pen}_n(f) $ that is at least
\begin{align*} 
\inf_{\widetilde{f}\in\widetilde{\mathcal{F}}} \{ \gamma_n L_n(\widetilde{f}) + n[P_n(T\widetilde{f}||f^{\star})-P_n(f||f^{\star}) + D^{\prime}_n(Tf,T\widetilde{f})] \}
\end{align*}
satisfies the risk bound
\begin{equation*}
\mathbb{E}\|T\hat{f}-f^{\star}\|^2 \leq (\tau+1)\inf_{f\in\mathcal{F}} \{ \|f-f^{\star}\|^2 + \mathbb{E}\text{pen}_n(f)/n + \mathbb{E}A_f \}.
\end{equation*}

By bounding the distortion in this way, we eliminate some error in approximating $ f $ by $ \widetilde{f} $ that arises from analyzing  $ P_n(T\widetilde{f}||f^{\star})-P_n(f||f^{\star}) $ and $ D_n(T\widetilde{f},f^{\star})- D_n(Tf,f^{\star}) $.

\begin{theorem} \label{thm:risk_rough}
Suppose $ \widetilde{\mathcal{F}} $ is a countable collection of functions that satisfies
\begin{equation*}
\mathbb{E}\sup_{\widetilde{f}\in\widetilde{\mathcal{F}}}\left\{ D_n^{\prime}(T\widetilde{f},f^{\star}) - \tau P_n(\widetilde{f}||f^{\star}) - \tau\gamma_n L_n(\widetilde{f}) \right\} \leq 0.
\end{equation*}
If $ \text{pen}_n(f) $ is at least
\begin{equation*} 
\inf_{\widetilde{f}\in\widetilde{\mathcal{F}}} \{ \gamma_n L_n(\widetilde{f}) + n[P_n(T\widetilde{f}||f^{\star})-P_n(f||f^{\star}) + D^{\prime}_n(Tf,T\widetilde{f})] \},
\end{equation*}
then the truncated estimator $ T\hat{f} $ with $ \hat{f} $ satisfying \prettyref{eq:estimator} has the resolvability bound
\begin{equation*} 
\mathbb{E}\|T\hat{f}-f^{\star}\|^2 \leq (\tau+1)\inf_{f\in\mathcal{F}} \{ \|f-f^{\star}\|^2 + \mathbb{E}\text{pen}_n(f)/n + \mathbb{E}A_f \}.
\end{equation*}
\end{theorem}

The main task is to construct the countable collection $ \widetilde{\mathcal{F}} $ and find a suitable upper bound on 
\begin{equation*} 
\inf_{\widetilde{f}\in\widetilde{\mathcal{F}}} \{ \gamma_n L_n(\widetilde{f}) + n[P_n(T\widetilde{f}||f^{\star})-P_n(f||f^{\star}) + D^{\prime}_n(Tf,T\widetilde{f})] \}.
\end{equation*}

\begin{lemma}
For every $ \widetilde{g} $, $ \widetilde{f} $, and $ f $,
\begin{align*}
& P_n(T\widetilde{f}||f^{\star})-P_n(f||f^{\star}) + D^{\prime}_n(Tf,T\widetilde{f}) \\
& \leq P_n(\widetilde{g}||f^{\star})-P_n(f||f^{\star}) + D^{\prime}_n(\widetilde{g},f) + 4\sqrt{2}B_n \sqrt{D(\widetilde{g}, \widetilde{f})+D^{\prime}_n(\widetilde{g},\widetilde{f})} + T_n.
\end{align*}
\end{lemma}

\begin{proof}
By \prettyref{lmm:truncate} \prettyref{eq:I} and \prettyref{eq:II},
\begin{align*}
(y-T\widetilde{f}_m(x))^2 - (y-f(x))^2
& = [(y-f_m(x))^2 - (y-f(x))^2] + \\ &
\qquad [(y-T\widetilde{f}_m(x))^2 -(y-Tf_m(x))^2] + \\ &
\qquad\qquad [(y-Tf_m(x))^2 - (y-f_m(x))^2] \\
& \leq [(y-f_m(x))^2 - (y-f(x))^2] + \\ &
\qquad 4B_n|f_m(x)-\widetilde{f}_m(x)| + \\ &
\qquad\qquad 4B_n(|y|-B_n)\mathbb{I}\{ |y| > B_n \} + \\ &
\qquad\qquad\qquad 2(|y|-B_n)^2\mathbb{I}\{ |y| > B_n \} \\
& = [(y-f_m(x))^2 - (y-f(x))^2] + \\ &
\qquad 4B_n|f_m(x)-\widetilde{f}_m(x)| + \\ &
\qquad\qquad 2(|y|^2- B_n^2)\mathbb{I}\{ |y| > B_n \}.
\end{align*}

By \prettyref{lmm:truncate} \prettyref{eq:III},
\begin{align*}
(T\widetilde{f}_m(x^{\prime})-Tf(x^{\prime}))^2
& \leq (f(x^{\prime})-f_m(x^{\prime}))^2 + 4B_n|\widetilde{f}_m(x^{\prime})-f_m(x^{\prime})|.
\end{align*}

Thus we find that $ (y-T\widetilde{f}_m(x))^2 - (y-f(x))^2 + (T\widetilde{f}_m(x^{\prime})-Tf(x^{\prime}))^2 $ is not greater than
\begin{align*}
[(y-f_m(x))^2 - (y-f(x))^2]+(f(x^{\prime})-f_m(x^{\prime}))^2 + & \\ \qquad 4B_n[|f_m(x)-\widetilde{f}_m(x)|+|\widetilde{f}_m(x^{\prime})-f_m(x^{\prime})|] + 2(|y|^2- B_n^2)\mathbb{I}\{ |y| > B_n \}
\end{align*}
\end{proof}

Recall that $ g $ is equal to $ f - f^{\star} $. In this way, there is a one to one correspondence between $ f $ and $ g $. To simplify notation, we sometimes write $ D_n(f,f^{\star}) $ as $ D_n(g) $ and $ D^{\prime}_n(f,f^{\star}) $ as $ D^{\prime}_n(g) $. Moreover, assume an analogous notation holds for the relative loss functions $ P_n(f||f^{\star}) $ and $ P^{\prime}(f||f^{\star}) $ and complexities $ L_n(f) $.
\begin{theorem} \label{thm:fundamental}
If $ \mathcal{F} $ is a countable collection of functions bounded in magnitude by $ B_n $ and $ L_n(f) $ satisfies the Kraft inequality $ \sum_{f\in\mathcal{F}}e^{-L_n(f)} \leq 1 $, then
\begin{equation*}
\mathbb{E}\sup_{f\in\mathcal{F}}\left\{ D_n^{\prime}(f||f^{\star}) - \tau P_n(f||f^{\star}) - \tau\gamma_n L_n(f)/n \right\} \leq 0,
\end{equation*}
where $ \tau = (1+\delta_1)(1+\delta_2) $ and $ \gamma_n = (2\tau)^{-1}(1+\delta_1/2)(1+2/\delta_1)(B+B_n)^2 + 2(1+1/\delta_2)\sigma^2 + 2(B+B_n)\eta $.
\end{theorem}
\begin{proof}
Let $ s^2(g) $ be as in \prettyref{lmm:concentration1}. Since $ g^2 $ is non-negative, $ s^2(g) \leq D^{\prime}_n(g^2) + D_n(g^2) $. Moreover, since $ |f| \leq B_n $ and $ |f^{\star}| \leq B $, it follows that $ s^2(g) \leq (B+B_n)^2(D^{\prime}_n(g) + D_n(g)) $. Let $ \gamma_1 = A_1(B+B_n)^2/2 $ with $ A_1 $ to be specified later. By \prettyref{lmm:concentration1}, we have
\begin{align} \label{eq:bound1}
\mathbb{E}\sup_{g\in\mathcal{G}}\left\{ (1-1/A_1)D^{\prime}_n(g)- (1+1/A_1)D_n(g)-\frac{\gamma_1}{n}L(g) \right\} & \\ \leq
\mathbb{E}\sup_{g\in\mathcal{G}}\left\{ D^{\prime}_n(g)-D_n(g)-\frac{\gamma_1}{n}L(g)-\frac{1}{2\gamma_1} s^2(g)\right\} \leq 0
\end{align}
By \prettyref{lmm:concentration2}, we also know that
\begin{equation} \label{eq:bound2}
\mathbb{E}\sup_{g\in\mathcal{G}}\left\{ \frac{1}{n}\sum_{i=1}^n\varepsilon_ig(X_i) -\frac{\gamma_2}{n}L(g) -\frac{1}{A_2n}D_n(g) \right\} \leq 0,
\end{equation}
where $ \gamma_2 = A_2\sigma^2/2 + (B+B_n)\eta $.
Adding the expression in \prettyref{eq:bound1} to $ 2a > 0 $ times the expression in \prettyref{eq:bound2} and collecting terms, we find that $ 1+1/A_1+2a/A_2 $ should be equal to $ a $ in order for $ D_n(g) $ and $ \frac{1}{n}\sum_{i=1}^n\varepsilon_ig(X_i) $ to be added together to produce $ P_n(g) $. Thus we find that
\begin{equation*}
\mathbb{E}\sup_{g\in\mathcal{G}}\left\{ (1-1/A_1)D^{\prime}_n(g) - a(P_n(g) + \frac{\gamma_n}{n}L(g)) \right\} \leq 0,
\end{equation*}
where $ \gamma_n = \gamma_1/a + 2\gamma_2 $.
Choosing $ A_1 = 1 + 2/\delta_1 $, $ A_2 = 2(1+1/\delta_2) $, and $ \tau = (1+\delta_1)(1+\delta_2) $, we find that $ a = \tau(1-1/A_1) $. Dividing the resulting expression by $ 1-1/A_1 $ produces
\begin{equation*}
\mathbb{E}\sup_{g\in\mathcal{G}}\left\{ D^{\prime}_n(g) - \tau P_n(g) - \tau\gamma_n L(g)/n \right\} \leq 0.
\end{equation*}
\end{proof}

In general, the penalty should not depend on the unknown test data $ \underline{X}^{\prime} $. However if one seeks to describe the error of a fit $ \hat{f} $ trained with the data $ (\underline{X},\underline{Y}) $ at new data points $ \underline{X}^{\prime} $, a penalty that depends on $ \underline{X}^{\prime} $ is natural and fits in with the standard trans-inductive setting in machine learning \cite{Gammerman1998}. Since we have been assuming that the input design $ X $ is contained in a cube with side length at most one, the dependence on $ \underline{X}^{\prime} $ of the penalties as shown in the following lemmata can be ignored.

When we speak of empirical $ L^2 $ covers of $ \mathcal{H} $, we mean with respect to the empirical measure of $ \underline{X}\cup\underline{X}^{\prime} $ on both the training and test data. That is, empirical $ L^2 $ covers of $ \mathcal{H} $ are with respect to the squared norm $ [D(h,\widetilde{h})+D^{\prime}(h,\widetilde{h})]/2 $.

\begin{theorem} \label{thm:penalty}
Let $ f = \sum_h \beta_h h $. Let $ \widetilde{\mathcal{H}}_1 $ be an empirical $ L^2 $ $ \epsilon_1 $-net for $ \mathcal{H} $ of cardinality $ M_1 $. Let $ \widetilde{\mathcal{H}}_2 $ be an empirical $ L^2 $ $ \epsilon_2 $-net for $ \mathcal{H} $ of cardinality $ M_2 $. Suppose these empirical covers do not depend on the underlying data. There exists a subset $ \widetilde{\mathcal{F}} $ of $ \mathcal{F} $ with cardinality at most $ \binom{M_2+M_1+m_0}{M_1+m_0} $ such that for $ v \geq v_f $ and $ \widetilde{v} = v(1+M_1/m_0) $
\begin{align*}
P_n(T\widetilde{f}||f^{\star})- P_n(f||f^{\star}) + D^{\prime}_n(Tf,T\widetilde{f})
& \leq
\frac{2\widetilde{v}^2\epsilon^2_1}{m_0} + \frac{\widetilde{v}^2M_1}{2m^2_0} + 8B_n\widetilde{v}\epsilon_2 + \frac{T_n}{n},
\end{align*}
for some $ \widetilde{f} $ in $ \widetilde{\mathcal{F}} $.

There exists a subset $ \widetilde{\mathcal{F}} $ of $ \mathcal{F} $ with cardinality at most $ \binom{M_2+m_0}{m_0} $
such that

\begin{align*}
P_n(T\widetilde{f}||f^{\star})- P_n(f||f^{\star}) + D^{\prime}_n(Tf,T\widetilde{f})
& \leq
\frac{2vv_f}{m_0}  + 8B_nv\epsilon_2 + \frac{T_n}{n},
\end{align*}
and in the case of no noise with $ B_n \geq B $
\begin{align*}
D_n(\widetilde{f},f) + D^{\prime}_n(f,\widetilde{f})
& \leq
\frac{4vv_f}{m_0} + 4v^2\epsilon^2_2
\end{align*}
for some $ \widetilde{f} $ in $ \widetilde{\mathcal{F}} $.

If $ \phi $ satisfies \prettyref{ass:ass2}, then there exists a subset $ \widetilde{\mathcal{F}} $ of $ \mathcal{F} $ with cardinality at most $ \binom{2\binom{2d+m_0}{m_0}+m_1}{m_1} $ such that
\begin{align*}
P_n(T\widetilde{f}||f^{\star})- P_n(f||f^{\star}) + D^{\prime}_n(Tf,T\widetilde{f})
& \leq
\frac{2v v_f}{m_1} + \frac{L_2v^2_fv^2_0}{m_0} + \frac{L^2_2v^2_fv^4_0}{4m_0}
\\ & \qquad + \left(\frac{1}{n}\sum_{i=1}^n|Y_i|\right)\frac{v_fL_2 v^2_0}{m_0} + \frac{T_n}{n}.
\end{align*}
for some $ \widetilde{f} $ in $ \widetilde{\mathcal{F}} $.

\end{theorem}

\begin{proof}
The proof is an immediate consequence of \prettyref{lmm:inequality1} and \prettyref{lmm:inequality}.
\end{proof}

According to \prettyref{thm:risk_rough} and \prettyref{thm:fundamental}, a valid penalty is at least
\begin{equation*} 
\gamma_n L_n(\widetilde{f}) + n[P_n(T\widetilde{f}||f^{\star})-P_n(f||f^{\star}) + D^{\prime}_n(Tf,T\widetilde{f})],
\end{equation*}
where $ \widetilde{f} $ belongs to a countable set $ \widetilde{\mathcal{F}} $ satisfying $ \sum_{\widetilde{f}\in\widetilde{\mathcal{F}}}e^{-L_n(\widetilde{f})} \leq 1 $. The constant $ \gamma_n $ is as prescribed in \prettyref{thm:fundamental}. By \prettyref{thm:penalty}, there is a set $ \widetilde{\mathcal{F}} $ with cardinality at most $ \binom{M_2+M_1+m_0}{M_1+m_0} $ such that for all $ f $ with $ v_f \leq v $, there is a $ \widetilde{f} $ in $ \widetilde{\mathcal{F}} $ such that $ P_n(T\widetilde{f}||f^{\star})-P_n(f||f^{\star}) + D^{\prime}_n(Tf,T\widetilde{f}) $ is bounded by
\begin{equation*}
\frac{2\widetilde{v}^2\epsilon^2_1}{m_0} + \frac{\widetilde{v}^2M_1}{2m^2_0} + 8B_n\widetilde{v}\epsilon_2 + \frac{T_n}{n}.
\end{equation*}
Using the fact that the logarithm of $ \binom{M_2+M_1+m_0}{M_1+m_0} $ is bounded by $ (M_1+m_0)\log(e(M_2/M_1+1)) $, a valid penalty divided by sample size is at least
\begin{align} \label{eq:noise_general}
\frac{\gamma_n}{n}(M_1+m_0)\log(e(M_2/M_1+1)) + \frac{2\widetilde{v}^2\epsilon^2_1}{m_0} + \frac{\widetilde{v}^2M_1}{2m^2_0} + 8B_n\widetilde{v}\epsilon_2 + \frac{T_n}{n}.
\end{align}

Alternatively, there is a set $ \widetilde{\mathcal{F}} $ with cardinality at most $ \binom{M_2+m_0}{m_0} $ such that for all $ f $ with $ v_f \leq v $, there is a $ \widetilde{f} $ in $ \widetilde{\mathcal{F}} $ such that $ P_n(T\widetilde{f}||f^{\star})-P_n(f||f^{\star}) + D^{\prime}_n(Tf,T\widetilde{f}) $ is bounded by
\begin{equation*}
\frac{2vv_f}{m_0} + 8B_nv\epsilon_2 + \frac{T_n}{n}
\end{equation*}
and hence a valid penalty divided by sample size is at least
\begin{equation} \label{eq:noise}
\frac{\gamma_n m_0\log M_2}{n} + \frac{2vv_f}{m_0} + 8B_nv\epsilon_2 + \frac{T_n}{n}.
\end{equation}

In the no noise case, a valid penalty divided by sample size is at least
\begin{equation} \label{eq:no_noise}
\frac{\gamma_nm_0\log M_2}{n} + \frac{4vv_f}{m_0} + 4v^2\epsilon^2_2.
\end{equation}

Analogously, if $ \phi $ satisfies \prettyref{ass:ass2}, a valid penalty divided by sample size is at least
\begin{equation} \label{eq:pen_general}
\frac{5m_0 m_1 \log(d+1)}{n} + \frac{2v v_f}{m_1} + \frac{L_2v^2_fv^2_0}{m_0} + \frac{L^2_2v^2_fv^4_0}{4m^2_0}
+ \left(\frac{1}{n}\sum_{i=1}^n|Y_i|\right)\frac{v_fL_2 v^2_0}{m_0} + \frac{T_n}{n}.
\end{equation}
for some $ \widetilde{f} $ in $ \widetilde{\mathcal{F}} $.

We now discuss how $ m_0 $, $ m_1 $, $ \epsilon_1 $, and $ \epsilon_2 $ should be chosen to produce penalties that yield optimal risk properties for $ T\hat{f} $.

\section{Risk bounds in high dimensions}

\subsection{Noise case under \prettyref{ass:ass1}}

By \prettyref{lmm:covering}, an empirical $ L^2 $ $ \epsilon_2 $-cover of $ \mathcal{H} $ has cardinality less than $ \binom{2d+\lceil ({v_0}/\epsilon_2)^2 \rceil }{\lceil ({v_0}/\epsilon_2)^2 \rceil} $. The logarithm of $ \binom{2d+\lceil ({v_0}/\epsilon_2)^2 \rceil }{\lceil ({v_0}/\epsilon_2)^2 \rceil} $ is bounded by $ 4({v_0}/\epsilon_2)^2\log(d+1) $.

Continuing from the expression \prettyref{eq:noise}, we find that $ \text{pen}_n(f)/n $ is at least
\begin{equation*} 
\frac{4\gamma_n m_0({v_0}/\epsilon_2)^2\log(d+1)}{n} + \frac{2vv_f}{m_0} + 8B_nv\epsilon_2 + \frac{T_n}{n}.
\end{equation*}

Choosing $ m_0 $ to be the ceiling of $ \left(\frac{vv_fn\epsilon^2_2}{2\gamma_n{v^2_0}\log(d+1)}\right)^{1/2} $, we see that $ \text{pen}_n(f)/n $ must be at least
\begin{equation*}
\frac{8\gamma_n{v^2_0}\log(d+1)}{n\epsilon^2_2} + 8\left(\frac{vv_f\gamma_n{v^2_0}\log(d+1)}{n\epsilon^2_2}\right)^{1/2}+8B_nv\epsilon_2 + \frac{T_n}{n}.
\end{equation*}
Finally, we set $ v = v_f $ and $ \epsilon_2 = \left(\frac{\gamma_n{v^2_0}\log(d+1)}{n B_n^2}\right)^{1/4} $ so that $ \text{pen}_n(f)/n $ must be at least
\begin{equation*}
16v_f\left(\frac{\gamma_n B_n^2{v^2_0}\log(d+1)}{n}\right)^{1/4} + 8\left(\frac{\gamma_n B_n^2{v^2_0}\log(d+1)}{n}\right)^{1/2} + \frac{T_n}{n}.
\end{equation*}

We see that the main term in the penalty divided by sample size is 
\begin{equation*}
16v_f\left(\frac{\gamma_n B_n^2{v^2_0}\log(d+1)}{n}\right)^{1/4}.
\end{equation*}

\subsection{No noise case under \prettyref{ass:ass1}}

Continuing from the expression \prettyref{eq:no_noise}, we find that $ \text{pen}_n(f)/n $ is at least
\begin{equation*} 
\frac{4\gamma_nm_0({v_0}/\epsilon_2)^2\log(d+1)}{n} + \frac{4vv_f}{m_0} + 4v^2\epsilon^2_2.
\end{equation*}

Choosing $ m_0 $ to be the ceiling of $ \left(\frac{vv_fn\epsilon^2_(\tau+1)^2}{\gamma_n{v^2_0}\log(d+1)(\tau+2)}\right)^{1/2} $, we see that $ \text{pen}_n(f)/n $ must be at least
\begin{equation*}
\frac{4\gamma_n{v^2_0}\log(d+1)}{n\epsilon^2_2} + 8\left(\frac{vv_f\gamma_n{v^2_0}\log(d+1)}{n\epsilon^2_2}\right)^{1/2}+ 4v^2\epsilon^2_2.
\end{equation*}

Finally, we set $ v = v_f $ and $ \epsilon_2 = \left(\frac{4(\tau+2)/(\tau+1)^2\gamma_n{v^2_0}\log(d+1)}{nv^2}\right)^{1/6} $ so that $ \text{pen}_n(f)/n $ must be at least
\begin{equation*}
16v_f^{4/3}\left(\frac{\gamma_n{v^2_0}\log(d+1)}{n}\right)^{1/3}+ 4(v_f^{4/3}+1)\left(\frac{\gamma_n{v^2_0}\log(d+1)}{n}\right)^{2/3},
\end{equation*}
where we used the fact that $ v^{2/3} \leq v^{4/3}+1 $.
We see that the main term in the penalty divided by sample size is 
\begin{equation*} 16v_f^{4/3}\left(\frac{\gamma_n{v^2_0}\log(d+1)}{n}\right)^{1/3}.
\end{equation*}

\subsection{In general under \prettyref{ass:ass2}}

Looking at \prettyref{eq:pen_general} suggests that we choose $ m_0 $ to be the floor of $ v^2_0 m_1 $ which results in a penalty divided by sample size of at least
\begin{equation*}
\frac{5\gamma_n m_0 m_1 \log(d+1)}{n} + \frac{2v^2_f}{m_1} + \frac{L_2v^2_fv^2_0}{m_0} + \frac{L^2_2v^2_fv^4_0}{4m^2_0}
+ \left(\frac{1}{n}\sum_{i=1}^n|Y_i|\right)\frac{v_fL_2 v^2_0}{m_0} + \frac{T_n}{n},
\end{equation*}
with leading terms of order
\begin{equation*}
\frac{\gamma_n v^2_0  m^2_1 \log(d+1)}{n} + \frac{v^2_f}{m_1}.
\end{equation*}
Choosing $ m_1 $ to be the floor of $ \left(\frac{v^2_f n}{\gamma_n v^2_0 \log(d+1)}\right)^{1/3} $ yields the conclusion that a valid penalty divided by sample size is at least of order
\begin{equation*}
v_f^{4/3}\left(\frac{\gamma_n{v^2_0}\log(d+1)}{n}\right)^{1/3} + v_f\left(\frac{1}{n}\sum_{i=1}^n|Y_i|\right)\left(\frac{\gamma_n{v^2_0}\log(d+1)}{n}\right)^{1/3} +  \frac{T_n}{n}.
\end{equation*}

\section{Risk bounds with improved exponents for moderate dimensions}

Continuing from the expression \prettyref{eq:noise_general}, we find that $ \text{pen}_n(f)/n $ is at least
\begin{equation*}
\frac{\gamma_n}{n}(M_1+m_0)\log(e(M_2/M_1+1)) + \frac{2\widetilde{v}^2\epsilon^2_1}{m_0} + \frac{\widetilde{v}^2M_1}{2m^2_0} + 8B_n\widetilde{v}\epsilon_2 + \frac{T_n}{n}.
\end{equation*}
Note that we can bound $ B^2_n $ by $ 
\gamma_n $ by choosing $ \delta_1 $ and $ \delta_2 $ appropriately. For the precise definition of $ \gamma_n $, see \prettyref{thm:fundamental}.
The strategy for optimization is to first consider the terms
\begin{equation} \label{eq:portion1}
\frac{\gamma_n}{n}m_0\log(e(M_2/M_1+1)) + \frac{2\widetilde{v}^2\epsilon^2_1}{m_0} + 8\sqrt{\gamma_n}\widetilde{v}\epsilon_2.
\end{equation}
After $ m_0 $, $ M_1 $, and $ M_2 $ have been selected, we then check that
\begin{equation} \label{eq:portion2}
\frac{\gamma_n}{n}M_1\log(e(M_2/M_1+1))+\frac{\widetilde{v}^2M_1}{2m^2_0}
\end{equation}
is relatively negligible. Choosing $ m_0 $ to be the ceiling of $ \left(\frac{2\widetilde{v}^2n\epsilon^2_1}{\gamma_n\log(e(M_2/M_1+1))}\right)^{1/2} $, we see that \prettyref{eq:portion1} is at most
\begin{equation*}
\frac{\gamma_n}{n}\log(e(M_2/M_1+1)) + 4\left(\frac{\widetilde{v}^2\gamma_n\epsilon^2_1\log(e(M_2/M_1+1))}{n}\right)^{1/2} + 8\sqrt{\gamma_n}\widetilde{v}\epsilon_2.
\end{equation*}
Note that an empirical $ L^2 $ $ \epsilon $-cover of $ \mathcal{H} $ has cardinality between $ ({v_0}/\epsilon)^d $ and $ (2{v_0}/\epsilon+1)^d \leq (3{v_0}/\epsilon)^d  $ whenever $ \epsilon \leq {v_0} $. Thus $ M_2/M_1 \leq (3\epsilon_1/\epsilon_2)^d $ whenever $ \epsilon_2 \leq {v_0} $ and hence
\begin{equation*}
\log(e(M_2/M_1+1)) \leq 1+(d/2)\log(9\epsilon^2_1/\epsilon^2_2+1) \leq d\log(9\epsilon^2_1/\epsilon^2_2+1),
\end{equation*} whenever $ \epsilon^2_1 \geq \epsilon^2_2(e-1)/9 $. These inequalities imply that \prettyref{eq:portion1} is at most
\begin{equation*}
\frac{d\gamma_n\log(9\epsilon^2_1/\epsilon^2_2+1)}{n} + 4\left(\frac{\widetilde{v}^2\epsilon^2_1d\gamma_n\log(9\epsilon^2_1/\epsilon^2_2+1)}{n}\right)^{1/2} + 8\sqrt{\gamma_n}\widetilde{v}\epsilon_2.
\end{equation*}
Next, set
\begin{equation*}
\epsilon^2_2 = \frac{9d\epsilon^2_1}{n}.
\end{equation*}
This means that the assumption $ \epsilon^2_1 \geq \epsilon^2_2(e-1)/9 $ is valid provided $ d \leq n/(e-1) $.
 Thus $ \prettyref{eq:portion1} $ is at most
\begin{equation*}
\frac{d\gamma_n\log(n/d+1)}{n} + 20\epsilon_1\widetilde{v}\sqrt{\frac{d\gamma_n\log(n/d+1)}{n}}.
\end{equation*}

Next, we add in the terms from \prettyref{eq:portion2}. The selections of $ m_0 $ and $ \epsilon_1 $ make \prettyref{eq:portion2} at most
\begin{equation*}
\frac{M_1d\gamma_n\log(n/d+1)}{n}+\frac{M_1d\gamma_n\log(n/d+1)}{n\epsilon_1^2}
\end{equation*}
Since $ M_1 \leq (3{v_0}/\epsilon_1)^d $ whenever $ \epsilon_1 \leq {v_0} $, we find that \prettyref{eq:portion2} is at most

\begin{equation*}
\frac{(3{v_0})^dd\gamma_n\log(n/d+1)}{n\epsilon_1^{d}} + \frac{(3{v_0})^dd\gamma_n\log(n/d+1)}{n\epsilon_1^{d+2}}
\end{equation*}

Let $ \epsilon_1 = 3{v_0}\left(\frac{d\gamma_n\log(n/d+1)}{n}\right)^{1/2(d+3)} $. Choosing $ \widetilde{v} = v_f $, we see that a valid penalty divided by sample size is at least

\begin{align*}
60v_f{v_0}\left(\frac{d\gamma_n\log(n/d+1)}{n}\right)^{1/2+1/2(d+3)}+\frac{1}{{v^2_0}}\left(\frac{d\gamma_n\log(n/d+1)}{n}\right)^{1/2+1/2(d+3)} & \\
+ \left(\frac{d\gamma_n\log(n/d+1)}{n}\right)^{1/2+3/2(d+3)} + \frac{d\gamma_n\log(n/d+1)}{n} + \frac{T_n}{n}.
\end{align*}

Note that for the form of the above penalty to be valid, we need $ \frac{d\gamma_n\log(n/d)}{n} $ to be small enough to ensure that $ \epsilon_1 $ and $ \epsilon_2 $ are both less than $ {v_0} $.




\section{Proofs of the lemmata}

An important aspect of the above covers $ \widetilde{\mathcal{F}} $ is that they only depend on the data $ (\underline{X},\underline{X}^{\prime}) $ through $ \|\underline{X}\|^2_{\infty} + \|\underline{X}^{\prime}\|^2_{\infty} $, where $ \|\underline{X}\|^2_{\infty} = \frac{1}{n}\sum_{i=1}^n\|X_i\|^2_{\infty} $. Since the coordinates of $ \underline{X} $ and $ \underline{X}^{\prime} $ are restricted to belong to $ [-1,1]^d $, the penalties and quantities satisfying Kraft's inequality do not depend on $ \underline{X} $ and $ \underline{X}^{\prime} $. This is an important implication for the following empirical process theory. On the other hand, using the fact that $ \|\underline{X}\|^2_{\infty} + \|\underline{X}^{\prime}\|^2_{\infty} $ is symmetric in the coordinates of $ \underline{X} $ and $ \underline{X}^{\prime} $ and has a mean that is at most logarithmic in $ d $, the following bounds can be adapted to handle covers $ \widetilde{\mathcal{F}} $ that depend on the training and test data without imposing sup-norm controls.

\begin{lemma} \label{lmm:concentration1}
Let $ (\underline{X},\underline{X}^{\prime}) = (X_1,\dots,X_n,X_1^{\prime},\dots,X_n^{\prime}) $, where $ \underline{X}^{\prime} $ is an independent copy of the data $ \underline{X} $ and where $ (X_1,\dots,X_n) $ are component-wise independent but not necessarily identically distributed. A countable function class $ \mathcal{G} $ and complexities $ L(g) $ satisfying $ \sum_{g\in\mathcal{G}}e^{-L(g)} \leq 1 $ are given. Then for arbitrary positive $ \gamma $,
\begin{equation} \label{eq:expectation}
\mathbb{E}\sup_{g\in\mathcal{G}}\left\{ D^{\prime}_n(g)-D_n(g)-\frac{\gamma}{n}L(g)-\frac{1}{2\gamma} s^2(g)\right\} \leq 0,
\end{equation}
where $ s^2(g) = \frac{1}{n}\sum_{i=1}^n(g^2(X_i)-g^2(X^{\prime}_i))^2 $.
\end{lemma}
\begin{proof}
Let $ \underline{Z} = (Z_1,\dots,Z_n) $ be a sequence of independent centered Bernoulli random variables with success probability $ 1/2 $. Since $ X_i $ and $ X^{\prime}_i $ are identically distributed, $ g^2(X_i)-g^2(X^{\prime}_i) $ is a symmetric random variable and hence sign changes do not affect the expectation in \prettyref{eq:expectation}. Thus the right hand side of the inequality in \prettyref{eq:expectation} is equal to
\begin{equation*}
\mathbb{E}_{\underline{Z},\underline{X},\underline{X}^{\prime}}\sup_{g\in\mathcal{G}}\left\{\frac{1}{n}\sum_{i=1}^nZ_i(g^2(X_i)-g^2(X^{\prime}_i)) - \frac{\gamma}{n}L(g)-\frac{1}{2\gamma}s^2(g)\right\}.
\end{equation*}

Using the identity $ x = \lambda\log(x/\lambda) $ with $ \lambda = \gamma/n $, conditioning on $ \underline{X} $ and $ \underline{X}^{\prime} $, and applying Jensen's inequality to move $ \mathbb{E}_{\underline{Z}} $ inside the logarithm, we have that
\begin{align*}
\mathbb{E}_{\underline{Z}}\sup_{g\in\mathcal{G}}\left\{\frac{1}{n}\sum_{i=1}^nZ_i(g^2(X_i)-g^2(X^{\prime}_i)) - \frac{\gamma}{n}L(g)-\frac{1}{2\gamma}s^2(g)\right\} & \\ \leq
\frac{\gamma}{n}\log\mathbb{E}_{\underline{Z}}\sup_{g\in\mathcal{G}}\exp\left\{\frac{1}{\gamma}\sum_{i=1}^nZ_i(g^2(X_i)-g^2(X^{\prime}_i)) - L(g)-\frac{n}{2\gamma^2}s^2(g)\right\}.
\end{align*}
Replacing the supremum with the sum and using the linearity of expectation, the above expression is not more than
\begin{align*}
\frac{\gamma}{n}\log\sum_{g\in\mathcal{G}}\mathbb{E}_{\underline{Z}}\exp\left\{\frac{1}{\gamma}\sum_{i=1}^nZ_i(g^2(X_i)-g^2(X^{\prime}_i)) - L(g)-\frac{n}{2\gamma^2}s^2(g)\right\} & \\ =
\frac{\gamma}{n}\log\sum_{g\in\mathcal{G}}\exp\left\{ - L(g)-\frac{n}{2\gamma^2}s^2(g) \right\}\mathbb{E}_{\underline{Z}}\exp\left\{\frac{1}{\gamma}\sum_{i=1}^nZ_i(g^2(X_i)-g^2(X^{\prime}_i))\right\}.
\end{align*}
Next, note that by the independence of $ Z_1,\dots,Z_n $,
\begin{align*}
\mathbb{E}_{\underline{Z}}\exp\left\{\frac{1}{\gamma}\sum_{i=1}^nZ_i(g^2(X_i)-g^2(X^{\prime}_i))\right\}
& = \prod_{i=1}^n\mathbb{E}_{Z_i}\exp\left\{\frac{1}{\gamma}Z_i(g^2(X_i)-g^2(X^{\prime}_i))\right\}.
\end{align*}
Using the inequality $ e^x+e^{-x} \leq 2e^{x^2/2} $, each $ \mathbb{E}_{Z_i}\exp\left\{\frac{1}{\gamma}Z_i(g^2(X_i)-g^2(X^{\prime}_i))\right\} $ is not more than $ \exp\left\{\frac{1}{2\gamma^2}(g^2(X_i)-g^2(X^{\prime}_i))^2\right\} $. Whence 
\begin{equation*}
\mathbb{E}_{\underline{Z}}\exp\left\{\frac{1}{\gamma}\sum_{i=1}^nZ_i(g^2(X_i)-g^2(X^{\prime}_i))\right\} \leq \exp\left\{\frac{n}{2\gamma^2}s^2(g)\right\}.
\end{equation*}
The claim follows from the fact that $ \frac{\gamma}{n}\log \sum_{g\in\mathcal{G}}e^{-L(g)} \leq 0 $.
\end{proof}

\begin{lemma} \label{lmm:concentration2}
Let $ \underline{\varepsilon} = (\varepsilon_1,\dots,\varepsilon_n) $ be conditionally independent random variables given $ \{X_i\}_{i=1}^n $, with conditional mean zero, satisfying Bernstein's moment condition with parameter $ \eta > 0 $. A countable class $ \mathcal{G} $ and complexities $ L(g) $ satisfying 
\begin{equation*} 
\sum_{g\in\mathcal{G}}e^{-L(g)} \leq 1
\end{equation*} 
are given. Assume a bound $ K $, such that $ |g(x)| \leq K $ for all $ g $ in $ \mathcal{G} $. Then
\begin{equation*}
\mathbb{E}\sup_{g\in\mathcal{G}}\left\{ \frac{1}{n}\sum_{i=1}^n\varepsilon_ig(X_i) -\frac{\gamma}{n}L(g) -\frac{1}{An}\sum_{i=1}^ng^2(X_i) \right\} \leq 0.
\end{equation*}
where $ A $ is an arbitrary constant and $ \gamma = A\sigma^2/2+Kh $.
\end{lemma}
\begin{proof}
Using the identity $ x = \lambda\log(x/\lambda) $ with $ \lambda = \gamma/n $, conditioning on $ \underline{X} $, and applying Jensen's inequality to move $ \mathbb{E}_{\underline{\varepsilon}} $ inside the logarithm, we have that
\begin{align*}
\mathbb{E}_{\underline{\varepsilon}|\underline{X}}\sup_{g\in\mathcal{G}}\left\{\frac{1}{n}\sum_{i=1}^n\varepsilon_ig(X_i) - \frac{\gamma}{n}L(g)-\frac{1}{An}\sum_{i=1}^ng^2(X_i)\right\} & \\ \leq
\frac{\gamma}{n}\log\mathbb{E}_{\underline{\varepsilon}|\underline{X}}\sup_{g\in\mathcal{G}}\exp\left\{\frac{1}{\gamma}\sum_{i=1}^n\varepsilon_ig(X_i) - L(g)-\frac{1}{\gamma A}\sum_{i=1}^ng^2(X_i)\right\}.
\end{align*}
Replacing the supremum with the sum and using the linearity of expectation, the above expression is not more than
\begin{align*}
\frac{\gamma}{n}\log\sum_{g\in\mathcal{G}}\mathbb{E}_{\underline{\varepsilon}|\underline{X}}\exp\left\{\frac{1}{\gamma}\sum_{i=1}^n\varepsilon_ig(X_i) - L(g)-\frac{1}{\gamma A}\sum_{i=1}^ng^2(X_i)\right\} & \\ =
\frac{\gamma}{n}\log\sum_{g\in\mathcal{G}}\exp\left\{ - L(g)-\frac{1}{\gamma A}\sum_{i=1}^ng^2(X_i)\right\}\mathbb{E}_{\underline{\varepsilon}|\underline{X}}\exp\left\{\frac{1}{\gamma}\sum_{i=1}^n\varepsilon_ig(X_i)\right\}.
\end{align*}
Next, note that by the independence of $ \varepsilon_1,\dots,\varepsilon_n $ conditional on $ \underline{X} $,
\begin{align*}
\mathbb{E}_{\underline{\varepsilon}|\underline{X}}\exp\left\{\frac{1}{\gamma}\sum_{i=1}^n\varepsilon_ig(X_i)\right\}
& = \prod_{i=1}^n\mathbb{E}_{\varepsilon_i|X_i}\exp\left\{\frac{1}{\gamma}\varepsilon_ig(X_i)\right\}.
\end{align*}
By \prettyref{lmm:MGF}, each $ \mathbb{E}_{\varepsilon_i|X_i}\exp\left\{\frac{1}{\gamma}\varepsilon_ig(X_i)\right\} $ is not more than $ \exp\left\{ \frac{\sigma^2g^2(X_i)}{2\gamma^2(1-\eta K/\gamma)} \right\} $. Whence 
\begin{align*}
\mathbb{E}_{\underline{\varepsilon}|\underline{X}}\exp\left\{\frac{1}{\gamma}\sum_{i=1}^n\varepsilon_ig(X_i)\right\} 
& \leq \exp\left\{\frac{\sigma^2\sum_{i=1}^ng^2(X_i)}{2\gamma^2(1-\eta K/\gamma)}\right\} \\ & = 
\exp\left\{\frac{1}{\gamma A}\sum_{i=1}^ng^2(X_i)\right\},
\end{align*}
where the last line follows from the definition of $ \gamma $. The proof is finished after observing that $ \frac{\gamma}{n}\log \sum_{g\in\mathcal{G}}e^{-L(g)} \leq 0 $.
\end{proof}

\begin{lemma} \label{lmm:approx1}
For $ f = \sum_h \beta_h h $ and $ f_0 $ in $ \mathcal{F} $, there is a choice of $ h_1,\dots,h_m $ in $ \mathcal{H} $ with $ f_m = (v/m)\sum_{k=1}^m h_k $, $ v \geq v_f $ such that
\begin{equation*}
\|f_m-f_0\|^2 - \|f_0 - f\|^2 \leq \frac{vv_f}{m}.
\end{equation*}
Moreover, the same bound holds for any convex combination of $ \|f_m-f_0\|^2 - \|f_0 - f\|^2 $ and $ \rho^2(f_m,f) $, where $ \rho $ is a possibly different Hilbert space norm.
\end{lemma}

\begin{proof}
Let $ H $ be a random variable that equals $ hv $ with probability $ \beta_h/v $ and zero with probability $ 1-v_f/v $. Let $ H_1,\dots,H_m $ be a random sample from the distribution defining $ H $. Then $ \overline{H} = \frac{1}{m}\sum_{j=1}^mH_j $ has mean $ f $ and furthermore the mean of $ \|f_m-f_0\|^2 - \|f_0 - f\|^2 $ is the mean is $ \|f-\overline{H}\|^2 $. This quantity is seen to be bounded by $ vv_f/m $. As a consequence of the bound holding on average, there exists a realization of $ f_m $ of $ \overline{H} $ (having form $ (v/m)\sum_{k=1}^m h_k $) such that $ \|f_m-f_0\|^2 - \|f_0 - f\|^2 $ is also bounded by $ Vv_f/m $.
\end{proof}

The next lemma is an extension of a technique used in \cite{Makovoz1996} to improve the $ L^2 $ error of an $ m $-term approximation of a function in $ L_{1,\mathcal{H}} $. The idea is essentially stratified sampling with proportional allocation \cite{Neyman1934} used in survey sampling as a means of variance reduction. In the following, we use the notation $ \|\cdot\| $ to denote a generic Hilbert space norm.

\begin{lemma} \label{lmm:approx}
Let $ \widetilde{\mathcal{H}} $ be an $ L^2 $ $ \epsilon_1 $-net of $ \mathcal{H} $ with cardinality $ M_1 $. For $ f = \sum_h\beta_h h $ and $ f_0 $ in $ \mathcal{F} $, there is a choice of $ h_1,\dots,h_m $ in $ \mathcal{H} $ with $ f_m = (1/m_0)\sum_{k=1}^mb_kh_k $, $ m \leq m_0 + M_1 $ and $ \|b\|_1 \geq v_f $ such that
\begin{equation*} \|f_0-f_m\|^2-\|f_0-f\|^2 \leq \frac{ vv_f\epsilon^2_1}{m_0}.
\end{equation*}

Moreover, there is an equally weighted linear combination $ f_m = (v/m_0)\sum_{k=1}^m h_k $, $ v \geq v_f $, $ m \leq m_0 + M_1 $ such that
\begin{equation*} \|f_0-f_m\|^2-\|f_0-f\|^2 \leq \frac{v^2\epsilon^2_1(1+M_1/m_0)}{m_0} + \frac{ v^2M_1}{4m^2_0}.
\end{equation*}

The same bound holds for any convex combination of $ \|f_m-f_0\|^2 - \|f_0 - f\|^2 $ and $ \rho^2(f_m,f) $, where $ \rho $ is a possibly different Hilbert space norm.

\end{lemma}
\begin{proof}
Suppose the elements of $ \widetilde{\mathcal{H}} $ are $ \widetilde{h}_1,\dots,\widetilde{h}_{M_1} $. Consider the $ M_1 $ sets (or ``strata")
\begin{equation*}
\widetilde{\mathcal{H}}_j = \{h\in\mathcal{H}:\|h-\widetilde{h}_j\|^2 \leq \epsilon^2_1 \},
\end{equation*}
$ j = 1,\dots,M_1 $. By working instead with disjoint sets $ \widetilde{\mathcal{H}}_j\setminus \bigcup_{1 \leq i\leq j-1}\widetilde{\mathcal{H}}_i $, $ \widetilde{\mathcal{H}}_0 = \emptyset $, that are contained in $ \widetilde{\mathcal{H}}_j $ and whose union is $ \mathcal{H} $, we may assume that the $ \widetilde{\mathcal{H}}_j $ form a partition of $ \mathcal{H} $.
Let $ M = m_0 + M_1 $ and $ v_j = \sum_{h\in\widetilde{\mathcal{H}}_j} \beta_h $. 
To obtain the first conclusion, define a random variable $ H_j $ to equal $ hv_j $ with probability $ \beta_h/v_j $ for all $ h\in\widetilde{H}_j $. let $ H_{1,j},\dots,H_{n_j,j} $ be a random sample of size $ N_j = \left\lceil{\frac{v_jM}{V}}\right \rceil $, where $ V = \frac{vM}{m_0} $ and $ v \geq v_f $, from the distribution defining $ H_j $. Note that the $ N_j $ sum to at most $ M $. Define $ g_j = \sum_{h\in\widetilde{\mathcal{H}}_j}\beta_h h $ and $ \overline{f} = \sum_{j=1}^{M_1}\frac{1}{N_j}\sum_{k=1}^{N_j}H_{k,j} $. Note that the mean of $ \overline{f} $ is $ f $. This means the expectation of $ \|f_0-\overline{f}\|^2-\|f_0-f\|^2 $ is the expectation of $ \|f-\overline{f}\|^2 $, which is equal to $ \sum_{j=1}^{M_1}\mathbb{E}\|H_j-g_j\|^2/N_j $. Now $ \mathbb{E}\|H_j-g_j\|^2/N_j $ is further bounded by 
\begin{equation*}
(V/M)\sum_{h\in\widetilde{\mathcal{H}}_j}\beta_h\inf_{h_j}\|h-h_j\|^2 \leq (V/M)\sum_{h\in\widetilde{\mathcal{H}}_j}\beta_h\|h-\widetilde{h}_j\|^2 \leq \frac{ v_jv\epsilon^2_1}{m_0}.
\end{equation*} 
The above fact was established by noting that the mean of a real-valued random variable minimizes its average squared distance from any point $ h_j $. Summing over $ 1 \leq j \leq M_1 $ produces the claim. Since this bound holds on average, there exists a realization $ f_m $ of $ \overline{f} $ (having form $ (1/m_0)\sum_{k=1}^m b_kh_k $ with $ \|b\|_1 \geq v_f $) such that $ \|f_0-f_m\|^2-\|f_0-f\|^2 $ is also bounded by $ \frac{ vv_f\epsilon^2_1}{m_0} $. \\

For the second conclusion, we proceed in a similar fashion. Suppose $ n_j $ is a random variable that equals $ \left\lceil{\frac{v_jM}{V}}\right \rceil $ and $ \left\lfloor{\frac{v_jM}{V}}\right \rfloor $ with respective probabilities chosen to make its average equal to $ \frac{v_jM}{V} $. Furthermore, assume $ n_1,\dots,n_{M_1} $ are independent. Define $ V_j = \frac{V}{M}n_j $. Since $ V_j  \leq v_j + \frac{V}{M} $, the $ V_j $ sum to at most $ V $.
Let $ H_j $ be a random variable that equals $ hv_j $ with probability $ \beta_h/v_j $ for all $ h \in \widetilde{\mathcal{H}}_j $. For each $ j $ and conditional on $ n_j $, let $ H_{1,j},\dots,H_{n_j,j} $ be a random sample of size $ N_j = n_j + \mathbb{I}\{n_j =0\} $ from the distribution defining $ H_j $. Note that the $ N_j $ sum to at most $ M $. Define $ g_j = \sum_{h\in\widetilde{\mathcal{H}}_j}\beta_h h $ and $ \overline{f} = \sum_{j=1}^{M_1}\frac{1}{N_j}\sum_{k=1}^{N_j}H_{k,j} $. Note that the conditional mean of $ \overline{H} $ given $ N_1,\dots,N_{M_1} $ is $ g = \sum_{j=1}^{M_1}(V_j/v_j)g_j $ and hence the mean of $ \overline{f} $ is $ f $. This means the expectation of $ \|f_0-\overline{f}\|^2-\|f_0-f\|^2 $ is the expectation of $ \|f-\overline{f}\|^2 $, which is equal to $ \sum_{j=1}^{M_1}\mathbb{E}\|H_j-(V_j/v_j)g_j\|^2/N_j + \mathbb{E}\|f-g\|^2 $ by the law of total variance. Now $ \mathbb{E}\|H_j-(V_j/v_j)g_j\|^2/N_j $ is further bounded by 
\begin{equation*}
(V/M)^2(n_j/v_j)\sum_{h\in\widetilde{\mathcal{H}}_j}\beta_h\inf_{h_j}\|h-h_j\|^2 \leq \frac{ v^2M\epsilon^2_1}{m^2_0}.
\end{equation*} 
The above fact was established by noting that the mean of a real-valued random variable minimizes its average squared distance from any point $ h_j $. Next, note that by the independence of the coordinates of $ v_1,\dots,v_{M_1} $ and the fact that $ V_j $ has mean $ v_j $,
\begin{equation*}
\mathbb{E}\|f-g\|^2 = \mathbb{E}\|\sum_{j=1}^{M_1}(V_j/v_j-1)g_j\|^2 = (V/M)^2\sum_{j=1}^{M_1}(\|g_j\|^2/v^2_j)\mathbb{V}(n_j).
\end{equation*}
Finally, observe that $ \|g_j\|^2 \leq v^2_j $ and $ \mathbb{V}(n_j) \leq 1/4 $ (a random variable whose range is contained in an interval of length one has variance bounded by $ 1/4 $). This shows that $ \mathbb{E}\|f-g\|^2 \leq  \frac{v^2M_1}{4m^2_0} $.
Since this bound holds on average, there exists a realization $ f_m $ of $ \overline{f} $ (having form $ (v/m_0)\sum_{k=1}^m h_k $) such that $ \|f_0-f_m\|^2-\|f_0-f\|^2 $ is also bounded by $ \frac{ v^2\epsilon^2_1(1+M_1/m_0)}{m_0} + \frac{v^2M_1}{4m^2_0} $.
\end{proof}

\begin{lemma} \label{lmm:approx2}
There is a collection of functions $ \widetilde{\mathcal{F}} $ with cardinality at most $ \binom{2\binom{2d+m_0}{m_0}+m_1}{m_1} \lesssim d^{m_0m_1} $ such that for each $ f(x) = \sum_{h} \beta_h h(x) = \sum_{h} \beta_h \phi(\theta_h \cdot x) $, there exists $ \widetilde{f} $ in $ \widetilde{\mathcal{F}} $ such that for any $ v \geq v_f $,
\begin{equation*}
\| \widetilde{f} - f\|^2 \leq \frac{vv_f}{m_1} + \frac{L^2_2v_f^2v^4_0}{4m^2_0}.
\end{equation*}
and
\begin{align*}
\|g -\widetilde{f}\|^2 - \|g - f\|^2
& \leq \frac{vv_f}{m_1} + \frac{L_2v_f(\|g\|_1+v_f)v^2_0}{m_0},
\end{align*}
provided $ \phi $ satisfies \prettyref{ass:ass2}. If $ \phi $ satisfies \prettyref{ass:ass1}, then
\begin{equation*}
\| \widetilde{f} - f\|^2 \leq \frac{vv_f}{m_1} + \frac{L^2_1v_f^2v^2_0}{m_0}.
\end{equation*}
\end{lemma}

\begin{proof}
Define a joint probability distribution $ (\widetilde{\theta}_H, H) $ as follows. Let $ \mathbb{P}[\widetilde{\theta}_{ih} = e_i\text{sgn}(\theta_{ih}) | H=h] = \frac{|\theta_{ih}|}{v_0} $ and $ \mathbb{P}[\widetilde{\theta}_{ih} = 0 | H=h] = 1-\frac{\|\theta_h\|_1}{v_0} $ for $ i = 1,2,\dots, d $ and $ \mathbb{P}[ H = h ] = \frac{|\beta_h|}{v} $ and $ \mathbb{P}[ H = 0 ] = 1-\frac{v_f}{v} $ for all $ h \in \mathcal{H} $ and $ v \geq v_f $.

Take a random sample $ \underline{H} = (H_1, H_2, \dots, H_{m_1}) $ from the distribution defining $ H $. For each $ j $, let $ \underline{\widetilde{\theta}}_j = (\widetilde{\theta}_{1, H_j},  \widetilde{\theta}_{2, H_j}, \dots, \widetilde{\theta}_{m_0, H_j}) $, $ \underline{\widetilde{\theta}}'_j = (\widetilde{\theta}'_{1, H_j},  \widetilde{\theta}'_{2, H_j}, \dots, \widetilde{\theta}'_{m_0, H_j}) $, and $ \underline{\widetilde{\theta}} = ( \underline{\widetilde{\theta}}_1 , \dots, \underline{\widetilde{\theta}}_{m_0}) $ be a random sample from the distribution defining $ \widetilde{\theta}_{H_j} $ and define
\begin{equation} \label{eq:form}
\tilde{f}_{m_0,m_1}(x) = \frac{v}{m_1}\sum_{j=1}^{m_1} \text{sgn}(\beta_{H_j})\phi\left(\frac{v_0}{m_0}\sum_{k=1}^{m_0} \widetilde{\theta}_{k, H_j} \cdot x \right).
\end{equation}
An important observation is that the average of $ \tilde{f}_{m_0,m_1} $ is $ \mathbb{E}\tilde{f}_{m_0,m_1} $ and hence by a similar argument to \prettyref{lmm:approx1}, there exists a realization of $ \tilde{f}_{m_0, m_1} $ such that
\begin{equation} \label{eq:bound1}
\| \tilde{f}_{m_0,m_1} - \mathbb{E}\tilde{f}_{m_0,m_1}\|^2 \leq \frac{vv_f}{m_1}.
\end{equation}
Furthermore, by the bias-variance decomposition,
\begin{equation} \label{eq:bound3}
\mathbb{E} \| \tilde{f}_{m_0,m_1} - f\|^2 = \| \tilde{f}_{m_0,m_1} - \mathbb{E}\tilde{f}_{m_0,m_1}\|^2 + \| f - \mathbb{E}\tilde{f}_{m_0,m_1}\|^2.
\end{equation}
The second term may be bounded as follows. First, note that
\begin{equation*}
\mathbb{E}\tilde{f}_{m_0,m_1}(x) = \sum_h \beta_h \mathbb{E}\phi\left(\frac{v_0}{m_0}\sum_{k=1}^{m_0} \widetilde{\theta}_{k, H_j} \cdot x \right).
\end{equation*}
By \prettyref{ass:ass2}, we have the pointwise bound
\begin{align}
|f(x) - \mathbb{E}\tilde{f}_{m_0,m_1}(x)| & = | \sum_h \beta_h \phi(\theta_h \cdot x) - \sum_h \beta_h \mathbb{E}\phi\left(\frac{v_0}{m_0}\sum_{k=1}^{m_0} \widetilde{\theta}_{k, h} \cdot x \right) | \nonumber \\
& \leq L_2\sum_h |\beta_h| \mathbb{E}\left| \frac{v_0}{m_0}\sum_{k=1}^{m_0} \widetilde{\theta}_{k, h} \cdot x - \theta_h \cdot x \right|^2 \leq L_2\sum_h |\beta_h| \frac{v^2_0\|x\|_{\infty}}{2m_0} \nonumber \\
& \leq \frac{L_2v_fv^2_0\|x\|^2_{\infty}}{2m_0} \leq \frac{L_2v_fv^2_0}{2m_0}. \label{eq:bound2}
\end{align}
Combining \prettyref{eq:bound1} and \prettyref{eq:bound2}, we have shown that there exists a realization of $ \tilde{f}_{m_0, m_1} $ such that
\begin{equation*}
\| \tilde{f}_{m_0,m_1} - f\|^2 \leq \frac{vv_f}{m_1} + \frac{L^2_2v^2_fv^4_0}{4m^2_0}.
\end{equation*}
For the second statement, we also use the bias-variance decomposition to write
\begin{align*}
\mathbb{E}\| g - \tilde{f}_{m_0,m_1} \|^2 - \| g - f \|^2 & = \mathbb{E}\| \tilde{f}_{m_0,m_1} - \mathbb{E}\tilde{f}_{m_0,m_1}\|^2 + \| g - \mathbb{E}\tilde{f}_{m_0,m_1}\|^2 - \| g - f \|^2 \\
& = \mathbb{E}\| \tilde{f}_{m_0,m_1} - \mathbb{E}\tilde{f}_{m_0,m_1}\|^2 + \langle f - \mathbb{E}\tilde{f}_{m_0,m_1}, 2g - f - \mathbb{E}\tilde{f}_{m_0,m_1} \rangle.
\end{align*}
As before, the first term is less than $ \frac{vv_f}{m_1} $. By \prettyref{eq:bound2}, $ |f(x) - \mathbb{E}\tilde{f}_{m_0,m_1}(x)| \leq \frac{L_2v_fv^2_0}{2m_0} $, and combining this with the pointwise bounds $ |f| \leq v_f $ and $ |\mathbb{E}\tilde{f}_{m_0,m_1}| \leq v_f $, we have
\begin{equation*}
|\langle f - \mathbb{E}\tilde{f}_{m_0,m_1}, 2g - f - \mathbb{E}\tilde{f}_{m_0,m_1} \rangle| \leq \frac{L_2v_f(\|g\|_1+v_f)v^2_0}{m_0}.
\end{equation*}
If $ \phi $ satisfies \prettyref{ass:ass1}, we use \prettyref{eq:bound3} together with the pointwise bound
\begin{align*}
|f(x) - \mathbb{E}\tilde{f}_{m_0,m_1}(x)| & = | \sum_h \beta_h \phi(\theta_h \cdot x) - \sum_h \beta_h \mathbb{E}\phi\left(\frac{v_0}{m_0}\sum_{k=1}^{m_0} \widetilde{\theta}_{k, h} \cdot x \right) | \nonumber \\
& \leq L_1\sum_h |\beta_h| \mathbb{E}\left| \frac{v_0}{m_0}\sum_{k=1}^{m_0} \widetilde{\theta}_{k, h} \cdot x - \theta_h \cdot x \right| \leq L_1\sum_h |\beta_h| \frac{v_0\|x\|_{\infty}}{\sqrt{m_0}} \nonumber \\
& \leq \frac{L_1v_fv_0\|x\|_{\infty}}{\sqrt{m_0}} \leq \frac{L_1v_fv_0}{\sqrt{m_0}},
\end{align*}
which yields
\begin{equation*}
\mathbb{E}\| \widetilde{f} - f\|^2 \leq \frac{vv_f}{m_1} + \frac{L^2_1v_f^2v^2_0}{m_0}.
\end{equation*}

By two applications of \prettyref{lmm:cardinality} with $ m = m_1 $ and $ M = \binom{2d+m_0}{m_0} $, the number of functions having the form \prettyref{eq:form} is at most $ \binom{2\binom{2d+m_0}{m_0}+m_1}{m_1} $.
\end{proof}

\begin{lemma} \label{lmm:inequality1}
Suppose $ \phi $ satisfies \prettyref{ass:ass2}. Let $ \underline{y} = \{y_i\}_{i=1}^n $ be a sequence of real numbers and let $ \underline{x} = \{x_i\}_{i=1}^n $ and $ \underline{x}^{\prime} = \{x^{\prime}_i\}_{i=1}^n $ be sequences $ d $-dimensional vectors. For $ f = \sum_h\beta_h h $ in $ \mathcal{F} $, there is a choice of $ \theta_1, \dots, \theta_{m_0} $ in the set of $ d $ standard basis vectors for $ \mathbb{R}^d $ with $ \widetilde{f}_{m_0, m_1} = \frac{v}{m_1}\sum_{j=1}^{m_1} \phi(\frac{v_0}{m_0}\sum_{k=1}^{m_0} \theta_k \cdot x) $, $ v \geq v_f $, such that
\begin{align*}
& \frac{1}{n}\sum_{i=1}^n(y_i-T\widetilde{f}_{m_0, m_1}(x_i))^2-\frac{1}{n}\sum_{i=1}^n(y_i-f(x_i))^2+\frac{1}{n}\sum_{i=1}^n(T\widetilde{f}_{m_0, m_1}(x^{\prime}_i))^2-Tf(x^{\prime}_i))^2 \\
& \leq \frac{v v_f}{m_1} + \frac{L_2v^2_fv^2_0}{m_0} + \left(\frac{1}{n}\sum_{i=1}^n|y_i|\right)\frac{L_2v_fv^2_0}{m_0} + \frac{2}{n}\sum_{i=1}^n(|y_i|^2-B^2_n)\mathbb{I}\{ |y_i| > B_n \}.
\end{align*}
If $ \widetilde{\mathcal{F}} $ denotes the collection of functions of the form $ \widetilde{f}_{m_0, m_1} $, then $ \widetilde{\mathcal{F}} $ has cardinality at most  $ \binom{2\binom{2d+m_0}{m_0}+m_1}{m_1} $.
\end{lemma}

\begin{lemma} \label{lmm:inequality}
Let $ \underline{y} = \{y_i\}_{i=1}^n $ be a sequence of real numbers and let $ \underline{x} = \{x_i\}_{i=1}^n $ and $ \underline{x}^{\prime} = \{x^{\prime}_i\}_{i=1}^n $ be sequences $ d $-dimensional vectors. Let $ \widetilde{\mathcal{H}}_1 $ be an empirical $ L^2 $ $ \epsilon_1 $-net for $ \mathcal{H} $ with cardinality $ M_1 $ and $ \widetilde{\mathcal{H}}_2 $ be an empirical $ L^2 $ $ \epsilon_2 $-net for $ \mathcal{H} $ with cardinality $ M_2 $. For $ f = \sum_h\beta_h h $ in $ \mathcal{F} $, there is a choice of $ \widetilde{h}_1,\dots,\widetilde{h}_m $ in $ \widetilde{\mathcal{H}}_2 $ with $ \widetilde{f}_m = (v/m_0)\sum_{k=1}^m \widetilde{h}_k $, $ v \geq v_f $, and $ m \leq m_0 + M_1 $, such that

\begin{align*}
& \frac{1}{n}\sum_{i=1}^n(y_i-T\widetilde{f}_m(x_i))^2-\frac{1}{n}\sum_{i=1}^n(y_i-f(x_i))^2 + \frac{1}{n}\sum_{i=1}^n(T\widetilde{f}_m(x^{\prime}_i))^2-Tf(x^{\prime}_i))^2 \\
& \leq \frac{2v^2\epsilon^2_1(1+M_1/m_0)}{m_0} + \frac{v^2M_1}{2m^2_0} + \\ & \qquad 8B_nv\left(1+M_1/m_0\right)\epsilon_2
+ \frac{2}{n}\sum_{i=1}^n(|y_i|^2-B^2_n)\mathbb{I}\{ |y_i| > B_n \}.
\end{align*}
If $ \widetilde{\mathcal{F}} $ denotes the collection of functions of the form $ \widetilde{f}_m $, then $ \widetilde{\mathcal{F}} $ has cardinality at most $ \binom{M_2+M_1+m_0}{M_1+m_0} $.

Moreover, there is a choice of $ \widetilde{h}_1,\dots,\widetilde{h}_{m_0} $ in $ \widetilde{\mathcal{H}}_2 $ with $ \widetilde{f}_{m_0} = (v/m_0)\sum_{k=1}^{m_0} \widetilde{h}_k $, $ v \geq v_f $ such that 
\begin{align*}
& \frac{1}{n}\sum_{i=1}^n(y_i-T\widetilde{f}_{m_0}(x_i))^2-\frac{1}{n}\sum_{i=1}^n(y_i-f(x_i))^2 + \frac{1}{n}\sum_{i=1}^n(T\widetilde{f}_{m_0}(x^{\prime}_i)-Tf(x^{\prime}_i))^2
\\ & \leq \frac{2vv_f}{m_0} + 8B_nv\epsilon_2 + \frac{2}{n}\sum_{i=1}^n(|y_i|^2-B^2_n)\mathbb{I}\{ |y_i| > B_n \}.
\stepcounter{equation}\tag{\theequation}\label{eq:ineq1}
\end{align*}
and
\begin{align*}
\frac{1}{n}\sum_{i=1}^n(\widetilde{f}_{m_0}(x_i)-f(x_i))^2 + \frac{1}{n}\sum_{i=1}^n(\widetilde{f}_{m_0}(x^{\prime}_i)-f(x^{\prime}_i))^2
& \leq \frac{4vv_f}{m_0} + 4v^2\epsilon^2_2.
\stepcounter{equation}\tag{\theequation}\label{eq:ineq2}
\end{align*}
If $ \widetilde{\mathcal{F}} $ denotes the collection of functions of the form $ \widetilde{f}_{m_0} $, then $ \widetilde{\mathcal{F}} $ has cardinality at most $ \binom{M_2+m_0}{m_0} $.

\end{lemma}

\begin{proof}
We only prove the first claim of the lemma. Inequalities \prettyref{eq:ineq1} and \prettyref{eq:ineq2} follow from similar arguments and \prettyref{lmm:approx1}.
Let $ f_m = (v/m_0)\sum_{k=1}^m h_k $ be as in the second part of \prettyref{lmm:approx}. Since $ \widetilde{\mathcal{H}}_2 $ is an empirical $ L^2 $ $ \epsilon_2 $-net for $ \mathcal{H} $, for each $ h_k $ there is an $ \widetilde{h}_k $ in $ \widetilde{\mathcal{H}}_2 $ such that 

\begin{equation*} 
\frac{1}{2n}\sum_{i=1}^n|h_k(x_i)-\widetilde{h}_k(x_i)|^2+\frac{1}{2n}\sum_{i=1}^n|h_k(x^{\prime}_i)-\widetilde{h}_k(x^{\prime}_i)|^2 \leq \epsilon^2_2.
\end{equation*}

Let $ \widetilde{f}_m = (v/m_0)\sum_{k=1}^m \widetilde{h}_k $.
By \prettyref{lmm:truncate} \prettyref{eq:I} and \prettyref{eq:II},
\begin{align*}
(y-T\widetilde{f}_m(x))^2 - (y-f(x))^2
& = [(y-f_m(x))^2 - (y-f(x))^2] + \\ &
\qquad [(y-T\widetilde{f}_m(x))^2 -(y-Tf_m(x))^2] + \\ &
\qquad\qquad [(y-Tf_m(x))^2 - (y-f_m(x))^2] \\
& \leq [(y-f_m(x))^2 - (y-f(x))^2] + \\ &
\qquad 4B_n|f_m(x)-\widetilde{f}_m(x)| + \\ &
\qquad\qquad 4B_n(|y|-B_n)\mathbb{I}\{ |y| > B_n \} + \\ &
\qquad\qquad\qquad 2(|y|-B_n)^2\mathbb{I}\{ |y| > B_n \} \\
& = [(y-f_m(x))^2 - (y-f(x))^2] + \\ &
\qquad 4B_n|f_m(x)-\widetilde{f}_m(x)| + \\ &
\qquad\qquad 2(|y|^2- B_n^2)\mathbb{I}\{ |y| > B_n \}.
\end{align*}

By \prettyref{lmm:truncate} \prettyref{eq:III},
\begin{align*}
(T\widetilde{f}_m(x^{\prime})-Tf(x^{\prime}))^2
& \leq (f(x^{\prime})-f_m(x^{\prime}))^2 + 4B_n|\widetilde{f}_m(x^{\prime})-f_m(x^{\prime})|.
\end{align*}

Thus we find that $ (y-T\widetilde{f}_m(x))^2 - (y-f(x))^2 + (T\widetilde{f}_m(x^{\prime})-Tf(x^{\prime}))^2 $ is not greater than
\begin{align*}
[(y-f_m(x))^2 - (y-f(x))^2]+(f(x^{\prime})-f_m(x^{\prime}))^2 + & \\ \qquad 4B_n[|f_m(x)-\widetilde{f}_m(x)|+|\widetilde{f}_m(x^{\prime})-f_m(x^{\prime})|] + 2(|y|^2- B_n^2)\mathbb{I}\{ |y| > B_n \}
\end{align*}

By the second conclusion in \prettyref{lmm:approx}, 
\begin{align*}
\frac{1}{n}\sum_{i=1}^n(y_i-f_m(x_i))^2 - \frac{1}{n}\sum_{i=1}^n(y_i-f(x_i))^2 + & \\  \frac{1}{n}\sum_{i=1}^n(f_m(x^{\prime}_i)-f(x^{\prime}_i))^2 \leq \frac{2v^2\epsilon^2_1(1+M_1/m_0)}{m_0} + \frac{v^2M_1}{2m^2_0}.
\end{align*} 
By the concavity of the square root function,
\begin{equation*}  \frac{1}{2n}\sum_{i=1}^n|h_k(x_i)-\widetilde{h}_k(x_i)|+\frac{1}{2n}\sum_{i=1}^n|h_k(x^{\prime}_i)-\widetilde{h}_k(x^{\prime}_i)|
\end{equation*}
is also no greater than $ \epsilon_2 $. Using this, we have that
\begin{equation*} 
\frac{1}{n}\sum_{i=1}^n|f_m(x_i)-\widetilde{f}_m(x_i)|+\frac{1}{n}\sum_{i=1}^n|f_m(x^{\prime}_i)-\widetilde{f}_m(x^{\prime}_i)|
\leq 2v\left(1+M_1/m_0\right)\epsilon_2.
\end{equation*}

The last conclusion about the cardinality of $ \widetilde{\mathcal{F}} $ follows from \prettyref{lmm:cardinality}.

\end{proof}

\begin{lemma} \label{lmm:covering}
Let $ \underline{x} = \{x_i\}_{i=1}^n $, where each $ x_i $ is a $ d $-dimensional vector in $ \mathbb{R}^d $. Define $ \|\underline{x}\|^2_{\infty} = \frac{1}{n}\sum_{i=1}^n\|x_i\|^2_{\infty} $. There is a subset $ \widetilde{\mathcal{H}} $ of $ \mathcal{H} $ with cardinality at most $ \binom{2d+m}{m} $ such that for each $ h(x) = \phi(x\cdot \theta) $ with $ \|\theta\|_1 \leq {v_0} $ in $ \mathcal{H} $, there is $ \widetilde{h}(x) = \phi(x\cdot \widetilde{\theta}) $ in $ \widetilde{\mathcal{H}} $ such that $ \frac{1}{n}\sum_{i=1}^n|h(x_i)-\widetilde{h}(x_i)|^2 \leq {v_0}\|\theta\|_1\|\underline{x}\|^2_{\infty}/m $.
\end{lemma}

\begin{proof}
By the Lipschitz condition on $ \phi $, it is enough to prove the bound for $ \frac{1}{n}\sum_{i=1}^n|\theta\cdot x_i - \widetilde{\theta}\cdot x_i|^2 $.
Let $ v $ be a random vector that equals $ e_j\text{sgn}(\theta_j){v_0} $ with probability $ |\theta_j|/{v_0} $, $ j = 1,2,\dots,d $ and equals the zero vector with probability $ 1- \|\theta\|_1/{v_0} $. Let $ v_1, v_2,\dots, v_m $ be a random sample from the distribution defining $ v $. Note that the average of $ \overline{\theta} = \frac{1}{m}\sum_{j=1}^{m} v_j $ is $ \theta $ and hence the average of each $ |\overline{\theta}\cdot x_i - \theta\cdot x_i|^2 $ is the variance of $ v\cdot x_i $ divided by $ m $. Taking the expectation of the desired quantity, we have 
\begin{align*}
\mathbb{E}\frac{1}{n}\sum_{i=1}^n|\overline{\theta}\cdot x_i - \theta\cdot x_i|^2 
& = \frac{1}{n}\sum_{i=1}^n\mathbb{E}|\overline{\theta}\cdot x_i - \theta\cdot x_i|^2 \\
& \leq \frac{1}{n}\sum_{i=1}^n \frac{\mathbb{E}(v\cdot x_i)^2}{m} \\
& \leq \frac{1}{n}\sum_{i=1}^n \frac{\|x_i\|^2_{\infty}{v_0}\|\theta\|_1}{m} \\
& = {v_0}\|\theta\|_1\|\underline{x}\|^2_{\infty}/m.
\end{align*}
Since this bound holds on average, there must exist a realization $ \widetilde{\theta} $ of $ \overline{\theta} $ for which the inequality is also satisfied.
Consider the collection of all vectors of the form 
\begin{equation*}
({v_0}/m)\sum_{j=1}^{m} u_j,
\end{equation*}
where $ u_j $ is any of the $ 2d+1 $ signed standard basis vectors including the zero vector. This collection has cardinality bounded by the number of non-negative integer solutions $ q_1,q_2,\dots,q_{2d+1} $ to 
\begin{equation*}
q_1 + q_2 + \cdots + q_{2d+1} = m.
\end{equation*}
This number is $ \binom{2d+m}{m} $ with its logarithm is bounded by $ m\log(e(2d/m+1)) $ or $ 2m\log(d+1) $. An important aspect of the log cardinality of this empirical cover is that it is logarithmic (and not linear) in the dimension $ d $. This small dependence on $ d $ is what produces desirable risk bounds when $ d $ is significantly greater than the available sample size $ n $.
\end{proof}

\begin{lemma} \label{lmm:MGF}
Let $ Z $ have mean zero and variance $ \sigma^2 $. Moreover, suppose $ Z $ satisfies Bernstein's moment condition with parameter $ \eta > 0 $. Then
\begin{equation}
\mathbb{E}(e^{tZ}) \leq \exp\left\{ \frac{t^2\sigma^2/2}{1-\eta|t|} \right\}, \quad |t| < 1/\eta.
\end{equation}
\end{lemma}

\begin{lemma} \label{lmm:truncate}

Define $ Tf = \min\{ B_n, |f| \}\text{sgn}f $. Then
\begin{enumerate}[(I)]
\item \label{eq:I} $
(y-Tf)^2 \leq (y-f)^2 + 2(|y|-B_n)^2\mathbb{I}\{ |y| > B_n \} $,
\item \label{eq:II}
$ (y-Tf)^2 \leq (y-T\widetilde{f})^2 + 4B_n|f-\widetilde{f}| + 4B_n(|y|-B_n)\mathbb{I}\{ |y| > B_n \} $,
and
\item \label{eq:III}
$ (T\widetilde{f}-Tf)^2 \leq (f-f_1)^2 + 4B_n|f_1-\widetilde{f}| $.
\end{enumerate}
\end{lemma}
\begin{proof}
\prettyref{eq:I}
Since $ (y-Tf)^2 = (y-f)^2 + 2(f-Tf)(2y-f-Tf) $, the proof will be complete if we can show that 
\begin{equation*}
 (f-Tf)(2y-f-Tf) \leq (|y|-B_n)^2\mathbb{I}\{ |y| > B_n\}.
\end{equation*}
Note that if $ |f| \leq B_n $, the left hand size of the above expression is zero. Thus we may assume that $ |f| > B_n $, in which case $ f-Tf = \text{sgn}f(|f|-B_n) $. Thus
\begin{align*}
(f-Tf)(2y-f-Tf)
& = 2y\text{sgn}f(|f|-B_n) - (|f|-B_n)(|f|+B_n) \\
& \leq 2|y|(|f|-B_n) - (|f|-B_n)(|f|+B_n).
\end{align*}
If $ |y| \leq B_n $, the above expression is less than $ -(|f|-B_n)^2 \leq 0 $. Otherwise, it is a quadratic in $ |f| $ that attains its global maximum at $ |f| = |y| $. This yields a maximum value of $ (|y|-B_n)^2 $. \\

\prettyref{eq:II} For the second claim, note that
\begin{equation*}
(y-Tf)^2 = (y-T\widetilde{f})^2 + (T\widetilde{f}-Tf)(2y-T\widetilde{f}-Tf).
\end{equation*}
Hence, we are done if we can show that
\begin{equation*}
(T\widetilde{f}-Tf)(2y-T\widetilde{f}-Tf) \leq 4B_n|f-\widetilde{f}| + 4B_n(|y|-B_n)\mathbb{I}\{ |y| > B_n \}.
\end{equation*}
If $ |y| \leq B_n $, then 
\begin{align*}
(T\widetilde{f}-Tf)(2y-T\widetilde{f}-Tf) 
& \leq 4B_n|T\widetilde{f}-Tf| \\
& \leq 4B_n|\widetilde{f}-f|.
\end{align*}
If $ |y| > B_n $, then 
\begin{align*}
(T\widetilde{f}-Tf)(2y-T\widetilde{f}-Tf) 
& \leq 2|T\widetilde{f}-Tf||y| + 2B_n|T\widetilde{f}-Tf| \\
& = 2|T\widetilde{f}-Tf|(|y|-B_n) + 4B_n|T\widetilde{f}-Tf| \\
& \leq 4B_n(|y|-B_n) + 4B_n|\widetilde{f}-f|.
\end{align*}
\prettyref{eq:III} For the last claim, note that
\begin{align*}
(T\widetilde{f}-Tf)^2 
& =  (T\widetilde{f}-Tf_1)^2 + [2T\widetilde{f}-Tf_1-Tf](Tf_1-Tf) \\
& \leq (T\widetilde{f}-Tf_1)^2 + 4B_n|Tf_1-Tf| \\
& \leq (\widetilde{f}-f_1)^2 + 4B_n|f_1-f|
\end{align*}
\end{proof}

\begin{lemma} \label{lmm:truncate_bound}
Let $ Y = f^{\star}(X) + \varepsilon $ with $ |f^{\star}(X)| \leq B $. Suppose  

\begin{enumerate}[(I)]
\item \label{eq:SE} $ \mathbb{E}e^{|\varepsilon|/\nu}  < +\infty $ or
\item \label{eq:SG} $ \mathbb{E}e^{|\varepsilon|^2/\nu}  < +\infty $ 
\end{enumerate}
for some $ \nu > 0 $. Then $ \mathbb{E}[(Y^2- B_n^2)\mathbb{I}\{|Y|>B_n\}] $ is at most 
\begin{enumerate}[(I)]
\item $ (4\nu^2/n)\mathbb{E}e^{|\varepsilon|/\nu} $ provided $ B_n > \sqrt{2}(B+\nu\log n) $ or 
\item $ (2\nu/n)\mathbb{E}e^{|\varepsilon|^2/\nu} $ provided $ B_n > \sqrt{2}(B+\sqrt{\nu\log n}) $.
\end{enumerate}
\end{lemma}
\begin{proof}
Under assumption \prettyref{eq:SE},
\begin{align*}
\mathbb{P}(Y^2- B_n^2 > t ) 
& = \mathbb{P}(|Y| > \sqrt{t+ B_n^2}) \\
& \leq \mathbb{P}(|\varepsilon| > \sqrt{t+ B_n^2} - B) \\
& \leq \mathbb{P}(|\varepsilon| > (1/\sqrt{2})(\sqrt{t}+B_n) - B) \\
& \leq e^{-\frac{1}{\nu}\sqrt{\frac{t}{2}}}e^{-\frac{1}{\nu }(\frac{B_n}{\sqrt{2}}-B)}\mathbb{E}e^{|\varepsilon|/\nu}.
\end{align*}
The last inequality follows from a simple application of Markov's inequality after exponentiation. Integrating the previous expression from $ t = 0 $ to $ t = +\infty $ ($ \int_{0}^{\infty}e^{-\frac{1}{\nu}\sqrt{\frac{t}{2}}}dt = 4\nu^2 $) yields an upper bound on $ \mathbb{E}[(Y^2- B_n^2)\mathbb{I}\{|Y|>B_n\}] $ that is at most $ (4\nu^2/n)\mathbb{E}e^{|\varepsilon|/\nu} $ provided $ B_n > \sqrt{2}(B+\nu\log n) $.

Under assumption \prettyref{eq:SG},
\begin{align*}
\mathbb{P}(Y^2- B_n^2 > t ) 
& = \mathbb{P}(|Y|^2 > t+ B_n^2) \\
& \leq \mathbb{P}(|\varepsilon|^2 > (1/2)(t+ B_n^2) - B^2) \\
& \leq e^{-\frac{t}{2\nu}}e^{-\frac{1}{v}(\frac{ B_n^2}{2}-B^2)}\mathbb{E}e^{|\varepsilon|^2/\nu}.
\end{align*}
The last inequality follows from a simple application of Markov's inequality after exponentiation. Integrating the previous expression from $ t = 0 $ to $ t = +\infty $ ($ \int_{0}^{\infty}e^{-\frac{t}{2\nu}}dt = 2\nu $) yields an upper bound on $ \mathbb{E}[(Y^2- B_n^2)\mathbb{I}\{|Y|>B_n\}] $ that is at most $ (2\nu/n)\mathbb{E}e^{|\varepsilon|^2/\nu} $ provided $ B_n > \sqrt{2}(B + \sqrt{\nu \log n}) \geq \sqrt{2(B^2+\nu\log n)} $.
\end{proof}

\begin{lemma} \label{lmm:cardinality}
The number of functions having the form $ \frac{v}{m}\sum_{k=1}^m h_k $, where $ h_k $ belong to a library of size $ M $ is at most $ \binom{M-1+m}{m} \leq \binom{M+m}{m}  $ with its logarithm bounded by $ m\log(e(M/m+1)) $.
\end{lemma}
\begin{proof}
Suppose the elements in the library are indexed from $ 1 $ to $ M $. Let $ q_i $ be the number of terms in $ \sum_{k=1}^m h_k $ of type $ i $. Hence the number of function of the form $ \frac{v}{m}\sum_{k=1}^m h_k $ is at most the number of non-negative integer solutions $ q_1,q_2,\dots,q_M $ to $ q_1 + q_2 + \cdots + q_M = m $. This number is $ \binom{M-1+m}{m} $ with its logarithm bounded by the minimum of $ m\log(e((M-1)/m+1)) $ and $ m\log M $.
\end{proof}

\begin{theorem} \label{thm:splineapprox}
Let $ f^{\star}(x) = \int_{\mathbb{R}^d}e^{ix\cdot\omega}\widetilde{f}(\omega)d\omega $ with $ v_{f^{\star},s} = \int_{\mathbb{R}^d}\|\omega\|^s_1|\widetilde{f}(\omega)|d\omega $ be an arbitrary target function. If $ v_{f^{\star},2} $ is  finite, there exists a linear combination of ridge ramp functions $ f_m(x) = \sum_{k=1}^m c_k (a_k \cdot x + b_k)_{+} $ with $ \|a_k\|=1 $ and $ |b_k| \leq 1 $ such that
\begin{equation*}
\|f^{\star} - x\cdot\nabla f^{\star}(0) - f^{\star}(0)-f_m\|^2 \leq \frac{16v^2_{f^{\star},2}}{m}.
\end{equation*}
Furthermore, if $ v_{f^{\star},3} $ is  finite, there exists a linear combination of squared ridge ramp functions $ f_m(x) = \sum_{k=1}^m c_k (a_k \cdot x + b_k)^2_{+} $ with $ \|a_k\|=1 $ and $ |b_k| \leq 1 $ such that
\begin{equation*}
\|f^{\star} - x^T H_{f^{\star}}(0) x/2 - x\cdot\nabla f^{\star}(0) - f^{\star}(0)-f_m\|^2 \leq \frac{16v^2_{f^{\star},3}}{m},
\end{equation*}
where $ H_{f^{\star}}(0) $ is the Hessian of $ f^{\star} $ at the point zero
\end{theorem}

\begin{proof}
If $ f^{\star} $ can be extended to a function on $ L^2(\mathbb{R}^d) $ with Fourier transform $ \widetilde{f} $, the function $ f^{\star}(x) - x\cdot\nabla f^{\star}(0) - f^{\star}(0) $ can be written as the real part of
\begin{equation} \label{eq:Fourier_rep}
\int_{\mathbb{R}^d}(e^{i\omega\cdot x}-i\omega\cdot x-1)\widetilde{f}(\omega)d\omega.
\end{equation}
If $ |z| \leq c $, we note the identity
\begin{equation*}
-\int_{0}^{c}[(z-u)_{+}e^{iu}+(-z-u)_{+}e^{-iu}]du = e^{iz}-iz-1.
\end{equation*}
If $ c = \|\omega\|_1 $, $ z = \omega\cdot x $, $ \alpha = \alpha(\omega) = \omega/\|\omega\|_1 $, and $ u = \|\omega\|_1t $, $ 0 \leq t \leq 1 $, we find that
\begin{align*}
-\|\omega\|^2_1\int_{0}^{1}[(\alpha\cdot x-t)_{+}e^{i\|\omega\|_1t}+(-\alpha\cdot x-t)_{+}e^{-i\|\omega\|_1t}]dt = & \\ e^{i\omega\cdot x}-i\omega\cdot x - 1.
\end{align*}
Multiplying the above by $ \widetilde{f}(\omega) = e^{ib(\omega)}|\widetilde{f}(\omega)| $, integrating over $ \mathbb{R}^d $, and applying Fubini's theorem yields
\begin{equation*}
f^{\star}(x) -x\cdot \nabla f^{\star}(0) - f^{\star}(0) = \int_{\mathbb{R}^d}\int_{0}^{1}g(t,\omega)dtd\omega,
\end{equation*} where 
\begin{align*} g(t,\omega) & = -[(\alpha\cdot x-t)_{+}\cos(\|\omega\|_1t+b(\omega)) + \\ &
\qquad (-\alpha\cdot x-t)_{+}\cos(\|\omega\|_1t-b(\omega))]\|\omega\|^2_1|\widetilde{f}(\omega)|.
\end{align*}
Consider a density on $ \{-1,1\}\times[0,1]\times\mathbb{R}^d $ defined by
\begin{equation*}
p(z,t,\omega) = |\cos(z\|\omega\|_1t+b(\omega))|\|\omega\|^2_1|\widetilde{f}(\omega)|/v
\end{equation*}
where
\begin{equation*}
v = \int_{\mathbb{R}^d}\int_{0}^1[|\cos(\|\omega\|_1t+b(\omega))|+|\cos(\|\omega\|_1t-b(\omega))|]\|\omega\|^2_1|\widetilde{f}(\omega)|dtd\omega \leq 2v_{f^{\star},2}.
\end{equation*}
Consider a random variable $ h(z,t,\alpha)(x) $ that equals
\begin{equation*}
(z\alpha\cdot x-t)_{+} \ s(zt,\omega),
\end{equation*}
where $ s(t,\omega) = -\text{sgn}\cos(\|\omega\|_1t+b(\omega)) $. Note that $ h(z, t,\alpha)(x) $ has the form $ \pm(\alpha\cdot x-t)_{+} $. We see that
\begin{align*}
f^{\star}(x) -x\cdot \nabla f^{\star}(0) - f^{\star}(0) =  & \\ v\int_{\{-1,1\}\times[0,1]\times\mathbb{R}^d}h(z, t,\alpha)(x)dp(z\times t\times \omega).
\end{align*}
One can obtain the final result by sampling $ (z_1,t_1,\omega_1),\dots,(z_m,t_m,\omega_m) $ randomly from $ p(z,t,\omega) $ and considering the average $ \frac{v}{m}\sum_{k=1}^mh(z_k,t_k,\omega_k) $. Note that since $ x = (x)_{+} - (-x)_{+} $, we can regard $ x\cdot\nabla f^{\star}(0) $ as belonging to the linear span of $ \{ x\mapsto z(\alpha\cdot x - t)_{+} : \|\alpha\|_1 = 1, \; 0 \leq t \leq 1, \; z \in \{-1,1\} \} $.
An easy argument shows that its variance is bounded by $ 16v^2_{f^{\star},2}/m $. This simple argument can be extended to higher order expansions of $ f^{\star} $. The function $ f^{\star}(x) - x^TH_{f^{\star}}(0)x/2 - x\cdot\nabla f^{\star}(0) - f^{\star}(0) $ can be written as the real part of
\begin{equation} \label{eq:Fourier_rep_2}
\int_{\mathbb{R}^d}(e^{i\omega\cdot x}+(\omega\cdot x)^2/2-i\omega\cdot x-1)\widetilde{f}(\omega)d\omega.
\end{equation}
As before, the integrand in \prettyref{eq:Fourier_rep_2} admits an integral representation by
\begin{equation*}
(i/2)\|\omega\|^3_1\int_{0}^{1}[(-\alpha\cdot x-t)^2_{+}e^{-i\|\omega\|_1t}-(\alpha\cdot x-t)^2_{+}e^{i\|\omega\|_1t}]dt.
\end{equation*}
Employing a sampling argument from an appropriately defined density, we are able to approximate $ f^{\star}(x) - x^TH_{f^{\star}}(0)x/2 - x\cdot\nabla f^{\star}(0) $ by a linear combinations of $ m $ second order spline functions (having bounded internal parameters) $ (\alpha\cdot x - t)^2_{+} $ with a squared error bounded by $ 16v^2_{f^{\star},3}/m $.
\end{proof}

\section*{Acknowledgements}
The risk bound framework of this paper is self-contained, yet borrows heavily from the unpublished manuscript \cite{Barron2008} written with Cong Huang and Gerald Cheang, which was judged to be too long for the Annals of Statistics. The Annals did show favorable inclination to receipt of a shorter version. We are very indebted to their contributions. The present paper with Jason M. Klusowski manages a compression of much of that work, corrected determination of improved exponents of the rate, and new improved rates for the high dimensional setting. The authors would also like to thank W. D. Brinda for useful discussions.

\bibliographystyle{plain}
\bibliography{references}

\newpage

\end{document}